\documentclass{amsart} 

\usepackage{amsfonts}
\usepackage{amsmath}
\usepackage{amsthm}
\usepackage{amssymb}
\usepackage[english]{babel}
\usepackage{graphics}
\usepackage{latexsym}
\usepackage{longtable} \usepackage{mathrsfs}
\usepackage{graphicx}

\newtheorem{theorem}{Theorem}
\newtheorem{corollary}[theorem]{Corollary}
\newtheorem{lemma}[theorem]{Lemma}
\newtheorem{proposition}[theorem]{Proposition}

\newtheorem{claim}[theorem]{Claim}
\newtheorem{example}[theorem]{Example}

\theoremstyle{definition}
\newtheorem{definition}[theorem]{Definition}
\newtheorem{remark}[theorem]{Remark}

 \newcommand{\lharpoonup}{-\!\!\!-\!\!\!\!\rightharpoonup}

\newcommand{\dist}{\textrm{dist}}

\newcommand{\BPL}{\medskip \noindent \textbf{Proof of Lemma} }
\newcommand{\BPC}{\medskip \noindent \textbf{Proof of Claim} }
\newcommand{\BPCO}{\medskip \noindent \textbf{Proof of Corollary} }
\newcommand{\BPP}{\medskip \noindent \textbf{Proof of Proposition} }
\newcommand{\BPT}{\medskip \noindent \textbf{Proof of Theorem} }

\newcommand{\sub}{\subseteq}
\newcommand{\set}{\setminus}
\newcommand{\by}{\times}

\newcommand{\mC}{\mathcal{C}}

\newcommand{\noi}{\noindent}
\newcommand{\ms}{\medskip}

\newcommand{\lmapsto}{\longmapsto}
\newcommand{\la}{\lambda}
\newcommand{\La}{\Lambda}
\newcommand{\e}{\varepsilon}
\newcommand{\de}{\delta}
\newcommand{\al}{\alpha}

\newcommand{\p}{\partial} \newcommand{\larrow}{\longrightarrow}

\newcommand{\R}{\mathbb{R}}
\newcommand{\B}{\mathbb{B}}
\newcommand{\Ca}{\mathcal{C}_{\alpha}} 
\newcommand{\N}{\mathbb{N}} 
\newcommand{\weight}{{(\R, e^{c \text{Id}})]^N}}
\newcommand{\eq}{U_{xx} - \nabla W ( U ) = - c\; U_x}
\newcommand{\action}{E_c (U) = \int_{\R}\left\{ \frac{1}{2}\big| U_x\big|^2 + W( U ) \right\}e^{cx} dx}

\newcommand{\bt}{\begin{theorem}}\newcommand{\et}{\end{theorem}}
\newcommand{\bl}{\begin{lemma}}\newcommand{\el}{\end{lemma}}
\newcommand{\be}{\begin{equation}}\newcommand{\ee}{\end{equation}}
\newcommand{\bc}{\begin{claim}}\newcommand{\ec}{\end{claim}}
\newcommand{\bp}{\begin{proof}}\newcommand{\ep}{\end{proof}}

\begin{document}

\title{Heteroclinic Travelling Waves of Gradient Diffusion Systems}


\author{\textsl{Nicholas D. Alikakos}}
\address{Department of Mathematics,
University of Athens Panepistimioupolis 11584, Greece \& IACM of
FORTH, Greece} \email{nalikako@math.uoa.gr}
\thanks{NDA was partially supported by
Kapodistrias Grant No.\ 70/4/5622 at the University of Athens.}

\author{\textsl{Nikolaos I. Katzourakis}}
\address{Department of Mathematics, University of Athens,
Panepistimioupolis 11584, Greece}
\email{nkatzourakis@math.uoa.gr}


\date{December 20, 2007 and, in revised form, November 30, 2008.}



\begin{abstract}
\noi We establish existence of travelling waves to the gradient
system $u_t = u_{zz} - \nabla W(u)$ connecting two minima of $W$
when $u : \R \times (0,\infty) \larrow \R^N$, that is, we establish
existence of a pair $(U,c) \in [C^2(\R)]^N \by (0,\infty)$,
satisfying
\[
\left\{\begin{array}{l}
  U_{xx} - \nabla W ( U ) = - c\; U_x\\
   U(\pm \infty) = a^{\pm},\\
\end{array}\right.
\]
where $a^{\pm}$ are local minima of the potential $W \in
C_{\textrm{loc}}^2(\R^N)$ with $W(a^-)< W(a^+)=0$ and $N \geq 1$.
Our method is variational and based on the minimization of the
functional $E_c (U) = \int_{\R}\Big\{ \frac{1}{2}|U_x|^2 + W( U )
\Big\}e^{cx} dx$ in the appropriate space setup. Following
Alikakos-Fusco \cite{A-F}, we introduce an artificial constraint to
restore compactness and force the desired asymptotic behavior, which
we later remove. We provide variational characterizations of the
travelling wave and the speed. In particular, we show that
$E_c(U)=0$.
\end{abstract}

\maketitle

\section{Introduction}

Assume we are given a potential $W \in C^2_{\textrm{loc}}(\R^N)$
with several local minima, in general at \emph{different levels}.
Let $a^+$, $a^-$ be local minima with $W(a^+)=0$, $W(a^-)<0$. We
consider the problem of existence of a solution $(U,c)$ to the
system
\begin{equation}
\label{problem} \left\{\begin{array}{l}
  \eq\\
   U(\pm \infty) = a^{\pm}\\
\end{array}\right.
\end{equation}
where $c>0$ and $U: \R \larrow \R^N$ is in $[C^2(\R)]^N$ connecting
$a^{\pm}$, the dimension being any $N \geq 1$. A typical potential
with two minima and $N=2$ is shown in Fig. 1. Solutions of problem
(\ref{problem}) are known as \emph{heteroclinic travelling waves}.
They are special solutions of the form $U(z-ct)=u(z,t)$ to the
diffusion system with gradient structure:
\begin{equation} \label{Diffusion system}
u_t = u_{zz} - \nabla W(u)\;,\;\;\;u=u(z,t) :\;\; \R \times
(0,\infty) \larrow \R^N,
\end{equation}
and in addition heteroclinic connections of the dynamical system
corresponding to the 2nd order ODE system $\eq$. Physically, problem
(\ref{problem}) can be interpreted as the Newtonian Law of motion
with force term $- \nabla (-W)$ due to the potential $-W$ and
dissipation (friction) term $-c U_x$. In this context, $U(x)$
represents the trajectory of an ideal unit mass particle going from
a global maximum to an other local maximum of $-W$, asymptotically
in time.
\[
\underset{\text{Fig.\ 1: Simulation of the standard 2-well $W$
deformed (exmpl.\ \ref{2D-Example}), having minima at different
levels}}{\includegraphics[scale=0.68]{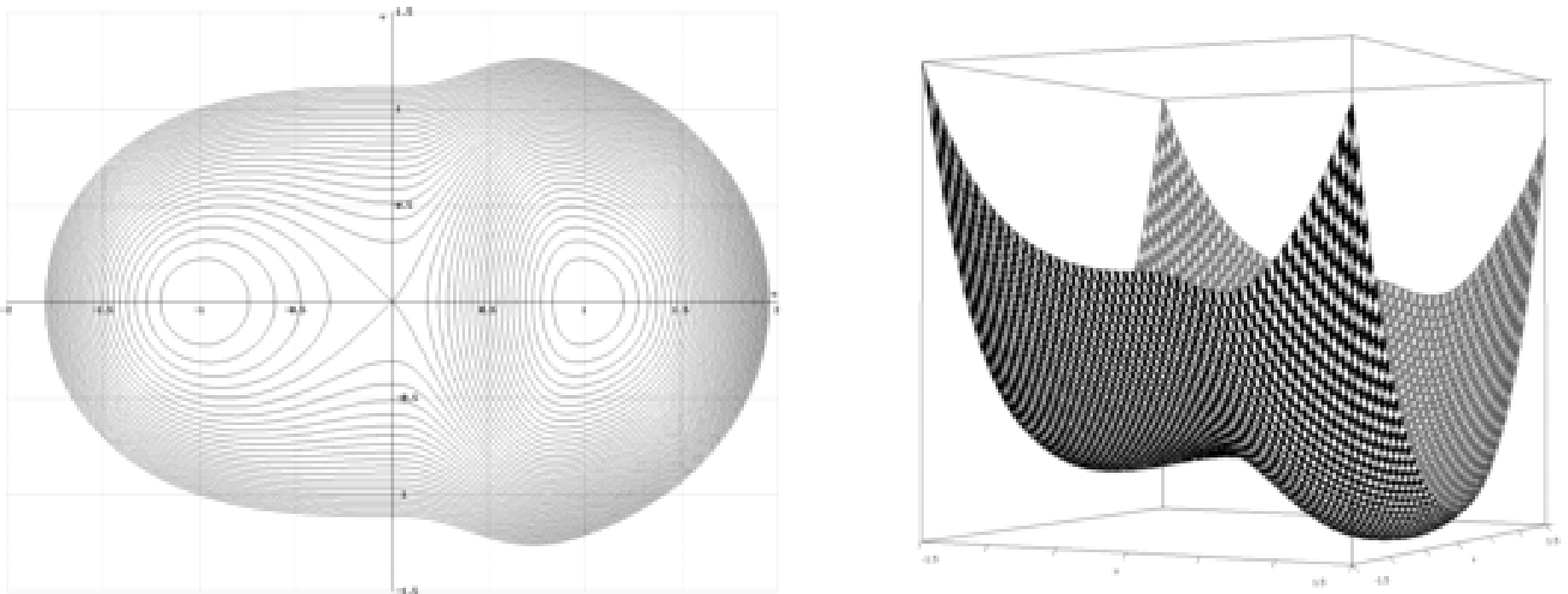}}
\]
Problem (\ref{problem}) with $c=0$ is a special case known as the
{\emph ``standing wave''} heteroclinic connection problem. It
reduces to a Hamiltonian system $U_{xx}=\nabla W(U)$ for a potential
with minima at the \emph{same level}. This case for general $N > 1$
has been studied by Sternberg in \cite{St}, Alikakos-Fusco in
\cite{A-F} and in great detail for $N=2$ by Alikakos, Betel\'u, Chen
in \cite{A-Be-C}.

\[
\underset{\text{Fig. 2: In general, no $\al$ - $\gamma$ connection
exists}}{\includegraphics[scale=0.14]{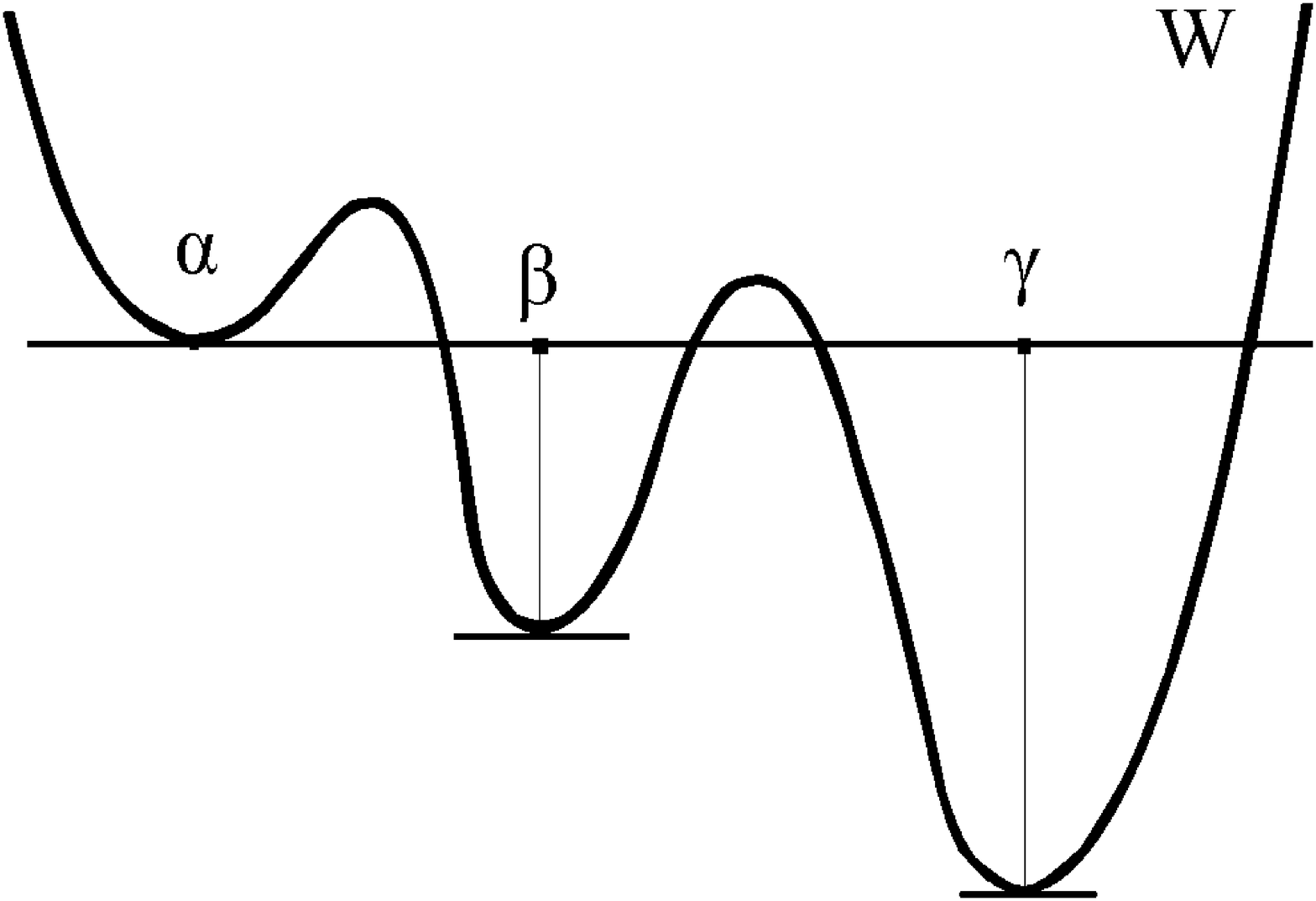}}
\]

The scalar case $N=1$ and $c>0$ of (\ref{problem}) is textbook
material from the viewpoint of existence (e.g \cite{He} p. 128,
\cite{Ev}, p. 175). The global stability of the connection for the
scalar case of (\ref{Diffusion system}) has been studied in the
classical papers of Fife and McLeod \cite{F-McL}, \cite{F-McL2} and
recently by Gallay and Risler in \cite{G-R}. Already in the scalar
case, existence for (\ref{problem}) of an heteroclinic between two
minima is not always guaranteed in the presence of a third one, as
it has been observed in \cite{F-McL} (Fig. 2).

In the vector case $N>1$ and $c \neq 0$ for (\ref{problem}) maximum
and comparison principles are no longer available and as a result
only special systems have been studied. We refer to the monograph of
Volperts' \cite{V} for monotone systems and numerous related
references.

In the very recent paper E. Risler \cite{R} has established
existence of solutions to (\ref{problem}), as a byproduct of his
study of the parabolic semiflow of (\ref{Diffusion system}). Among
other results, Risler studies the case of a bistable potential and
proves the existence of a travelling wave connecting the global
minimum of $W$ with a local minimum, as in the present paper.
However, his hypotheses are more restrictive than our \eqref{h**}
(Sec.\ 8), which shows the advantage of the Direct Method we
utilize.

Another very recent paper that establishes existence of travelling
waves, actually for a generalization of \eqref{problem} is
Lucia-Muratov-Novaga \cite{LMN}. Their method has similarities with
ours, but their hypotheses are different and not directly addressing
the potential $W$.

In the present paper we choose to work directly with the time
independent problem (\ref{problem}). We prove existence of
heteroclinic travelling waves for potentials with several minima
under weak coercivity requirements which allow for potentials
unbounded from below. We establish connections between possibly
degenerate minima, imposing assumptions only on the geometry of the
sublevel set $\big\{W \leq \al \big\}\subseteq \R^N$ for $\al>0$
small, which encloses the minima (assumptions \eqref{h*} in Sec.\ 6,
\eqref{h**} in Sec.\ 8).

Our approach is variational: we introduce a \emph{weighted action}
functional, an idea already introduced in Fife-McLeod (\cite{F-McL},
\cite{F-McL2}), to obtain travelling wave solutions to
(\ref{problem}) as (local) minimizers of the weighted action
\begin{equation} \label{action}
  \action
\end{equation}
in the \emph{Fr\'echet space} of vector functions
$[H^1_{\textrm{loc}}(\R)]^N$, utilizing certain devices to overcome
the unboundedness and compactness problems of $E_c$. We show that
action-minimizing travelling waves $(U,c)$ are characterized by the
property $E_c(U)=0$ and they can be derived as solutions to
\begin{equation} \label{Variational formulation}
E_c(U) = \inf\ \Big\{E_c(V) :\ V \in [H_{\textrm{loc}}^1(\R)]^N,\
V(\pm \infty )= a^{\pm} \Big\} \; , \;\;\;\;\;\;\ E_c(U) = 0 .
\end{equation}

We now give a brief description of our method. A formal computation
shows that critical points of $E_c$ correspond to weak solutions of
(\ref{problem}). We wish to construct solutions of $\eq$, with the
desired behavior at infinity $U(\pm \infty) = a^{\pm}$, by
minimizing (\ref{action}), in the appropriate setup. Minimization
can not be done directly, because the unbounded domain $\R$ excludes
strong compactness in all reasonable functional spaces, while the
asymptotic behavior required in (\ref{problem}) can not be
guaranteed.

In addition, (\ref{action}) is \emph{not} generally \emph{bounded}
from below for all $c>0$, a difficulty not present when $c=0$, and
moreover it is sensitive to translations: $E_c \big(U(\cdot - \de)
\big)= e^{c \de} E_c (U)$. Thus, a minimizing sequence may converge
to the trivial minimizers $a^{\pm}$ with $E_c(a^+)=0$,
$E_c(a^-)=-\infty$.

To overcome these problems, we first solve a constrained
minimization problem, utilizing the \emph{unilateral constraint
method} introduced by Alikakos and Fusco in \cite{A-F}: we fix 2
arbitrary parameters $c, L > 0$ and we minimize $E_c$ directly
within the \emph{admissible set} of functions in
$[H_{\textrm{loc}}^1(\R)]^N$ whose graph lies in the cylinders $(-
\infty , -L] \times \mathbb{B}(a^-, r_0) $ and $[L, + \infty) \times
\mathbb{B}(a^+, r_0) $ enclosing the 2 minima $a^{\pm}$ to be
connected. Minimization leads to a 2-parameter family of minimizers
in $c,L>0$. Then $L$ is increased with the hope that the constraint
is not realized for some minimizer, thus solving the Euler-Lagrange
equation (\ref{problem}) for some specific value of the other
parameter $c=c^*>0$.

This device bounds from below (\ref{action}), and allows us to
``capture'' an object which is close to a solution to
(\ref{problem}). Constrained minimizers are \emph{piecewise
solutions} (except possibly at the rims $\{\pm L\} \times
\partial\big( \mathbb{B}(a^{\pm},r_0)\big)$) converging
asymptotically to $a^{\pm}$, for all $c>0$. The main effort in the
proof is devoted to showing that the constraint is in fact
\emph{not} realized for a specific $c^*>0$ and for sufficiently
large $L$.

The role of "$c$" is as follows. We incorporate into $E_c$ an
\emph{arbitrary parameter} $c>0$ which, until Sec. 6, is always
\emph{arbitrary and fixed}. In particular, we do \emph{not} view $c$
as a functional $c(U)$ of $U$. The specific $c=c^*$ which guarantees
existence is determined by the requirement that $E_{c^*}(U_L)=0$ for
sufficiently large $L \geq L^*$. This is necessary for existence of
minimizers since translation sensitivity of (\ref{action}) shows
that the only possible \emph{finite} infimum of (\ref{action}) is
zero. A more transparent characterization was pointed out by the
referee and is as follows. First look for the smallest possible
value $c>0$ for which \eqref{action} is bounded form below over $\{U
\in [H^1_{\textrm{loc}}(\R)]^N : U(\pm \infty)=a^\pm \}$. Then, for
that $c$ construct the travelling wave by minimizing (\ref{action}).
A nice consequence of this is a uniqueness property of the speed for
minimizing travelling waves.

The paper is organized as follows. In  Sec.\ 2 we solve the
constraint problem for $E_c$ in $[H^1_{\textrm{loc}}(\R)]^N$,
resulting to a 2-parameter family of minimizers in $c>0$ and $L>0$.
In Sec.\ 3, assuming a very mild local monotonicity (h) near the
minima $a^{\pm}$, we show that constrained minimizers are piecewise
solutions to $\eq$, solving it on $\R \set \{\pm L\}$ and converging
to $a^{\pm}$ at ${\pm \infty}$.

In Sec.\ 4 we introduce the main tool for removing the constraint,
two local replacement lemmas, modeled after Lemmas 3.3, 3.4 in
\cite{A-F}. The new ingredient is the introduction of a convex set
in the place of a ball, which allows controlling the solution far
from the minima. The presentation here is self-contained independent
of the rest of the paper.

In Sec.\ 5 we establish certain energy identities. In particular,
they imply an energy equipartition at $+ \infty$ and that $E_c(U_L)$
measures the jumps $[[(U_{L})_x]]\big|_{\pm L}$.

In Sec.\ 6 we introduce a global assumption \eqref{h*} and determine
the speed $c^*$ of the travelling wave. $c^*$ is defined by means of
a variational formula (see \eqref{SET}) which is similar to a
formula of Heinze \cite{Hei}. Utilizing tools from Sec.\ 4, 5, we
prove that $c^*$ satisfies the desired properties (Proposition
\ref{Determination of $c=c^*$}). Hence, we distinguish the
\emph{suitable} $E_{c^*}$ among all $\big\{E_{c} : c>0 \big\}$. The
variational formulation (\ref{Variational formulation}) which
implies existence for (\ref{problem}) is also given here.

In Sec.\ 7 we prove existence of solution by removing the constraint
and derive explicit bounds on $c^* \in [c_{\min}, c_{\max}]$, by
means of our variational formulation (\ref{Variational
formulation}).

In Sec.\ 8 we show that the assumption \eqref{h*} can be relaxed to
include potentials that are unbounded from below or have other
critical points besides $a^\pm$ (cf. \cite{A-F}). Finally, in the
Appendix we discuss the optimality of our assumptions.

Our proof includes the special case $W(a^-)=W(a^+)=0$, $c=0$ that
was treated in \cite{A-F}.

\section{The Constrained Minimization Problem}

Here we solve a minimization problem for $\action$ in the local
Sobolev space $[H^1_{\textrm{loc}}(\R)]^N$ of vector $U : \R \larrow
\R^N$. $[H^1_{\textrm{loc}}(\R)]^N$ admits a Fr\'echet topology,
defined by the seminorms of $[H^1(-m,m)]^N$, $m\geq 1$. Technically,
instead of $[H^1_{\textrm{loc}}(\R)]^N$ we use its isomorphic copy
$[H^1_{\textrm{loc}}\weight$ with weight $x \mapsto e^{cx}$, the
standard Lebesgue measure $dx$ being replaced by the absolutely
continuous $e^{cx}dx$. It is only a matter of convenience, since
minimization gives derivatives bounded in $[L^2(\R, e^{c
\text{Id}})]^N$. $C^k_{\textrm{loc}}(\R^N)$ will denote the space of
$C^k$ functions equipped with the Fr\'echet topology of uniform
convergence together with all the derivatives over compacts, while
$C^k(\R^N)$ denotes the bounded $C^k$ functions with its standard
norm. We shall frequently decompose $W$ as $W = W^+ - W^-$, where
$W^+ =\max\{W,0\}$ and $W^- =\max \{-W,0\}$.

\bl \textbf{(Characterization of the speed)} \label{Characterization
of c} Assume that a solution $(U,c)$ to (\ref{problem}) exists,
satisfying $U_x(\pm \infty)=0$ up to sequences. Then:
\[
W^-(a^-) \ =\ c \int_{\R} \big|U_x \big|^2 dx
\;\;\;\;\;\;\&\;\;\;\;\; c\big(a^+ - a^-\big)\ =\  \int_{\R}\nabla
W\big(U\big)dx.
\]
\el
  \BPL \ref{Characterization of c}. The equation readily implies $-U_{xx} \cdot U_x \;+\; \nabla W(U) \cdot
U_x = c \big|U_x \big|^2$. Hence,
\begin{align*}
c \int _{\R} \big|U_x \big|^2dx \ & =\ - \int _{\R}\big(
\frac{1}{2}\big|U_x \big|^2\big)_x dx \;+\; \int_{\R} \big( W(U)
\big)_x dx\\
& =\  \pm \;0 \;+\ W(U(+ \infty)) \;-\; W(U(- \infty))\\
& = -\; W(a^-).
\end{align*}
Moreover, again from the equation we have
\begin{align*}
\int_{\R}\nabla W\big(U\big)dx \; &=\; \int_{\R}\Big(U_{xx} \;+\;c
U_x \Big)dx \\
& =\; 0 -\;0 \;+\; c\big(a^+ - a^-\big). \qed
\end{align*}
\noindent As a consequence of Lemma \ref{Characterization of c}, if
$U(\pm \infty) = a^{\pm}$ and $W(a^+) =0>W(a^-)$, then $c$ must be
positive.

\noi Take now $L>0$ and $r_0 >0$ small, such that $W(u) \geq 0$ for
$|a^+ - u| \leq r_0$ and $W(u) < 0$ for $|a^- - u| \leq r_0$. We
introduce the constraint sets:
\begin{align*}
\mathcal{X}^+_L \ :=& \  \Big\{ U \in
[H^1_{\textrm{loc}}\weight\;:\;|U(x)
- a^+| \leq r_0,\;x \geq +L \Big\},\\
 \mathcal{X}^-_L \ :=& \ \Big\{ U \in
[H^1_{\textrm{loc}}\weight\;:\;|U(x) - a^-| \leq r_0,\;x \leq -L
\Big\},
\end{align*}
and set $\mathcal{X}_L := X^+_L \bigcap \mathcal{X}^-_L $. Pointwise
values make sense by means of the imbedding
$[H^1_{\textrm{loc}}\weight \hookrightarrow
[C^0_{\textrm{loc}}(\R)]^N$.

\bt \textbf{(Existence of Constrained Minimizers)}
\label{Constrained Problem} Let $W$ be a potential in
$C^2_{\textrm{loc}}(\R^N)$ and $a^{\pm}$ two of its local minima,
with $W(a^-)<0=W(a^+)$, and $a^-$ its global minimum. We assume that
$W^{-1}\big([W(a^-),0]\big)$ is compact in $\R^N$. If $L>0$, $c>0$
are fixed parameters, then the minimization problem
\[
E_c(U_L)\ = \ \underset{\mathcal{X}_L }{\inf}\ \big\{E_c \big\}
\]
has a solution $U_L$ in $\mathcal{X}_L \sub
[H^1_{\textrm{loc}}\weight$. \et

\noi The assumption on $W$ implies ${\liminf}_{|u|\rightarrow
\infty}\big[W(u)\big] \geq 0$. This will be relaxed in the sequel,
allowing for potentials with several local minima and possibly
unbounded negative values, by means of a localization. We denote the
minimizers of $E_c$ into $\mathcal{X}_L$ by $U_L$ instead of the
more accurate notation $U_{c,L}$, suppressing the dependence on the
parameter $c>0$ which (until Sec. 6) is always fixed.

\BPT \ref{Constrained Problem}. We first show that $\mathcal{X}_L
\neq \emptyset$ together with $-\infty < {\inf}_{\mathcal{X}_L }
\big\{E_c\big\} < \infty$. Since we are interested only in
increasing the parameter $L$ later, we restrict, as we can, our
attention to $L \geq 1$.

\noi \textbf{Claim.} \emph{There exists an affine function $U_{aff} \in
\mathcal{X}_L \bigcap [W^{1,\infty}_{\textrm{loc}}(\R)]^N$ such that
\[
- \infty\ <  \ -\frac{e^{cL} W^-(a^-)}{c} \ \leq \;
\underset{\mathcal{X}_L }{\inf}\ \big\{E_c\big\} \; \leq \
E_c(U_{aff}) \ < \ \infty.
\]}

  \BPC. Let $\chi_A$
denote the characteristic of $A \subseteq \R$. We set
\[
U_{aff}(x)\;:=\; a^- \chi_{(- \infty , -1)} \;+\; \left( \frac{1 -
x}{2} a^- + \frac{1 + x}{2} a^+\right)\ \chi_{[-1 , 1]} \;+\;a^+
\chi_{(1, \infty)}.
\]
\[
\underset{\text{Fig.3\ The device of constrained minimization which
restores compactness and
boundedness}}{\includegraphics[scale=0.18]{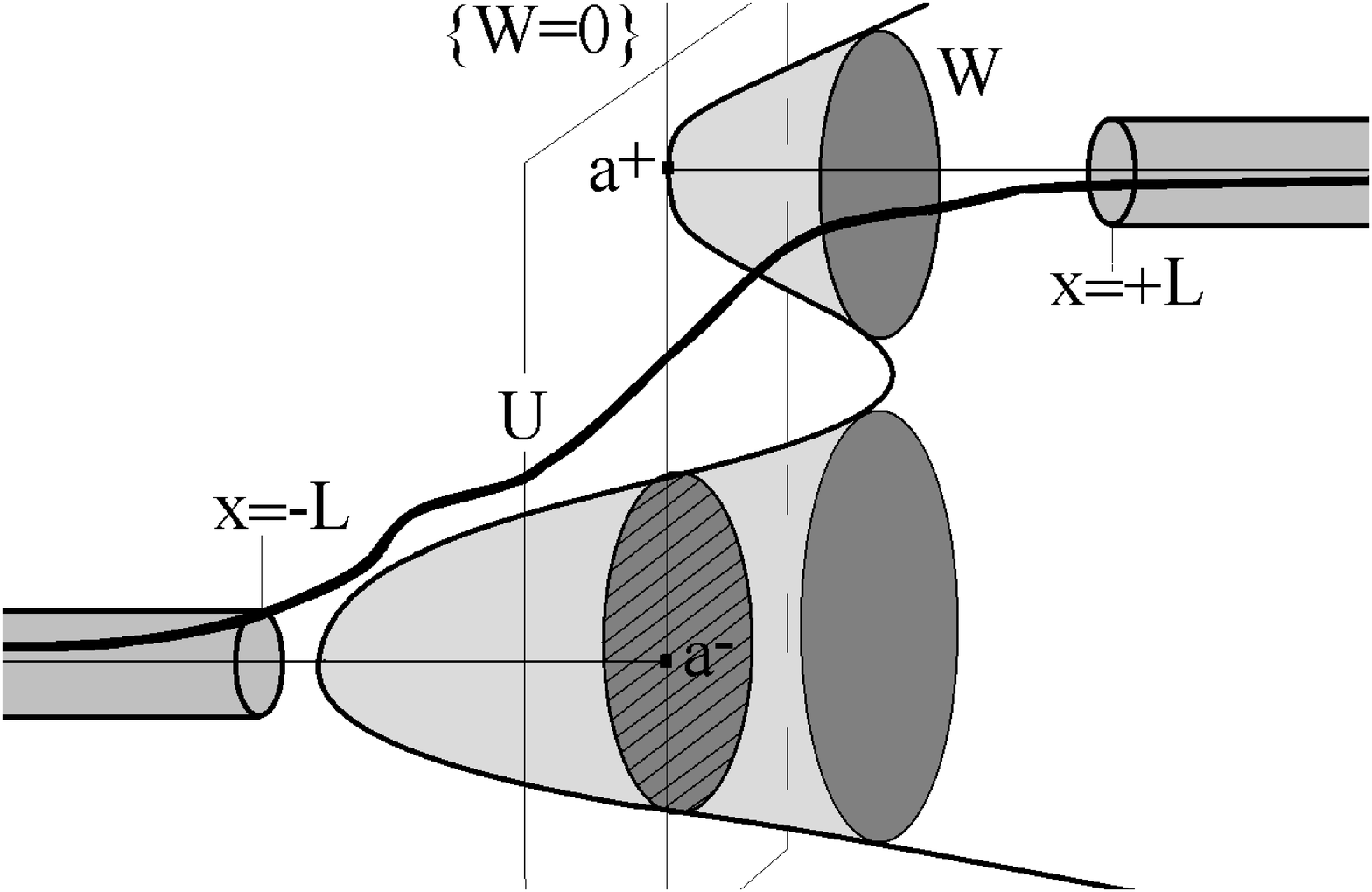}}
\]
Clearly, $(U_{aff})_x \in [L^\infty_{\textrm{loc}}(\R)]^N$ and
exists a.e. on $\R$. Thus, $U_{aff} \in [H^1_{\textrm{loc}}\weight$.
We calculate
\begin{align*}
E_c(U_{aff})\;=&\; \int_{- \infty}^{-1}(0+W(a^-))e^{cx}dx \;+\;
\int^{ \infty}_{1}(0+W(a^+))e^{cx}dx\\
 &+\ \int_{-1}^{1}\left\{
\frac{1}{2} \Big| \frac{a^+ - a^-}{2} \Big|^2 + W \left( \frac{1 -
x}{2} a^- + \frac{1 + x}{2} a^+ \right) \right\}e^{cx}dx\\
 \leq & \ \int_{-1}^{1}
 \left\{ \frac{1}{2} \Big| \frac{a^+ - a^-}{2} \Big|^2  \ +\  W^+ \left(
\frac{1 - x}{2} a^- + \frac{1 + x}{2} a^+ \right)
\right\} e^{cx}dx\\
 & \  + \ \frac{1}{c} e^{-c} W(a^-).
 \end{align*}
Hence, if we set $E^+_c(U):= \int_{\R}\Big\{ \frac{1}{2}| U_x|^2 +
W^+ \big( U \big) \Big\}e^{cx} dx$, we obtain
\begin{equation}
\label{Affine function bound} E_c(U_{aff}) \; \leq \; -e^{-c} \frac{
W^-(a^-)}{c} \;+\; e^{c} E^+_0(U_{aff}).
\end{equation}
This implies the upper bound $\sup_{L\geq 1}
{\inf}_{\mathcal{X}_L}\big\{E_c\big\} \leq \sup_{L\geq 1}
E_c(U_{aff}) < \infty$. If $U$ lies in $\mathcal{X}_L$, we have
$W^-(U(x))=0$ for $x\geq L$ and $W^+(U(x))=0$ for $x\leq -L$. Hence,
for any such $U$, utilizing that $W^-(U) \leq W^-(a^-)$, we have
\begin{align*}
E_c (U) &=\ \int_{\R}\left\{ \frac{1}{2}\big| U_x\big|^2 + W \big( U
\big) \right\}e^{cx} dx\\
 & = \  \frac{1}{2} \int_{\R}|U_x|^2
e^{cx}dx \;+\; \int_{\R} W^+ (U)e^{cx}dx \ -\;
\int_{\R}W^-(U)e^{cx}dx
\\
& \geq \;-\; \int_{\R} W^-(U)e^{cx}dx \\
&\geq \; - \ W^-(a^-) \int_{- \infty}^L e^{cx}dx \ = \; -
\frac{W^-(a^-)}{c}e^{cL} . \qed
\end{align*}

\noi By $C^2$ regularity of solutions to \eqref{problem}, we may
assume that ${\inf}\ {\mathcal{X}_L} \big[E_c\big] < E_c(U_{aff})$
strictly. We choose a minimizing sequence $\{U^n_L\}_{n \geq 1}$ in
$[H^1_{\textrm{loc}}\weight$ such that $E_c(U^n_L ) \larrow
{\inf}_{\mathcal{X}_L} \big\{E_c\big\}$, as $n \rightarrow \infty$.
The constraints immediately yield
\[
\big|U^n_L(x)\big|\;\leq\; \max\big\{|a^+|, |a^-|\big\} \;+\;
r_0\;,\;\;\;\;\;x \in \;(-\infty, -L] \cup [L, \infty).
\]
\noi \textbf{Claim.} \textbf{(Uniform Bounds)} \emph{There exists a
$C=C(c,L,W)>0$ such that
\[
\sup_{n\geq 1}\ \big\|(U^n_L)_x \big\|_{[L^2\weight}\ \leq \ C
\;,\;\;\; \sup_{n\geq 1}\ \big\|U^n_L \big\|_{[L^{\infty}(\R)]^N} \
\leq \ C.
\]}

\BPC. For any $x \in [-L,L]$, we have the estimates
\begin{align*}
|U^n_L(x)| \ &\leq\ |U^n_L(-L)| \;+\; \int_{-L}^x \big|(U^n_L)_t
\big| e^{\frac{ct}{2}} e^{\frac{-ct}{2}} dt \\
 &\leq \max\big\{|a^+|, |a^-|\big\} \;+\; r_0 \;+\;
\Bigg(\int_{-L} ^L e^{-ct}dt \Bigg)^{\frac{1}{2}}\Bigg( \int_{-L}^x
\big|(U^n_L)_t \big|^2 e^{ct} dt \Bigg)^{\frac{1}{2}},
\end{align*}
\begin{align*}
\frac{1}{2}\int_{\R} \big|(U^n_L)_x \big|^2 e^{cx}dx\ &\leq \
E_c(U_{aff})\;-\; \int_{\R} W(U^n_L) e^{cx}dx \\
& \leq\ E_c(U_{aff})- \int_{\R} W^+(U^n_L) e^{cx}dx \;+\; \int_{\R}
W^-(U^n_L) e^{cx}dx\\
& \leq \ E_c(U_{aff}) \;+\;  \int_{-\infty}^L W^-(U^n_L) e^{cx}dx\\
& \leq \ E_c(U_{aff}) \;+\; \frac{W^-(a^-)}{c}e^{cL}.
\end{align*}
We conclude:
\[
\frac{1}{2} \big\|(U^n_L)_x \big\|^2_{[L^2\weight} \;\leq\;
\frac{W^-(a^-)}{c}e^{cL} \;+\; E_c(U_{aff}).
\]
Utilizing that $\big|U^n_L(x)\big|\leq \max\big\{|a^+|, |a^-|\big\}
+ r_0$ for $x \in (-\infty, -L] \bigcup [L, \infty)$, we get
\[
\big\|U^n_L \big\|_{[L^{\infty}(\R)]^N} \;\leq \; \max\big\{|a^+|,
|a^-|\big\} \;+\; r_0 \;+\;
\Bigg(\frac{e^{cL}-e^{-cL}}{c}\Bigg)^{\frac{1}{2}}\big\|(U^n_L)_x
\big\|_{[L^2\weight}.  \qed
\]
We may now proceed to the existence of the minimizer. By the claim
above, $(U^n_L)_1^{\infty}$ is \emph{bounded in the locally convex
sense} in $[H^1_{\textrm{loc}}\weight$, with the derivatives bounded
in $[L^2(\R,e^{cId})]^N$:
\begin{align*}
{\sup}_{n \geq 1} \big\| U^n_L \big\|_{\big(H^1(I,e^{cId})]^N} & \
\leq \ C(c,L,W,I)\ \text{ for all $I
\subset \subset \R$, \ \ }\\
{\sup}_{n \geq 1} \big\| (U^n_L)_x \big\|_{[L^2(\R,e^{cId})]^N} & \
\leq \ C(c,L,W).
\end{align*}
By standard compactness arguments, there is a $U_L \in
[H^1_{\textrm{loc}}\weight$ such that up to a certain subsequence $
U^n_L \lharpoonup U_L$ as $n \rightarrow \infty$ weakly in
$[H^1_{\textrm{loc}}(\R)]^N$ and $U^n_L \larrow U_L$ in
$[L^2_{\textrm{loc}}\weight$ and a.e.\ on $\R$. By weak LSC of the
weighted $L^2$ norm and the Fatou Lemma for $W(U^n_L)+
W^-(a^-)\chi_{(-\infty, L]} \geq 0$\footnote{We owe this argument to
the referee.}, we have
\begin{align*}
  \int_{\R} \frac{1}{2} |(U_L)_x|^2  e^{cx}dx\
& \leq \ \underset{n \rightarrow \infty}{\liminf} \ \int_{\R}
\frac{1}{2} |(U^n_L)_x|^2  e^{cx}dx , \\
 \int_{\R} \Big\{ W(U_L)+  W^-(a^-) \chi_{(-\infty, L]} \Big\}
e^{cx}dx & \leq \\
\leq \  \underset{n \rightarrow \infty}{\liminf} & \ \int_{\R}
\Big\{W(U^n_L)+ W^-(a^-) \chi_{(-\infty, L]} \Big\}e^{cx}dx.
\end{align*}
Hence, the theorem follows together with the bounds
\begin{align*}
 -\frac{e^{cL} W^-(a^-)}{c} \;\leq \; E_c(U_L)
& \; \leq\; \underset{n \rightarrow \infty}{\lim} E_c(U^n_L) \\
& \;\leq \; -e^{-c} \frac{ W^-(a^-)}{c} \;+\; e^{c} E^+_0(U_{aff}).
\qed
\end{align*}

\section{Constrained Minimizers are Piecewise Solutions}

We will now prove that the constrained minimizers $U_L$ of Theorem
\ref{Constrained Problem} are piecewise solutions in
$C^2_{\textrm{loc}}(\R \set \{-L,L\})^N$, while $U_L|_{(- \infty,
-L]}$ and $U_L|_{[L, \infty)}$ are inside the cylinders and converge
asymptotically to $a^{\pm}$. Following \cite{A-F}, we introduce the
following local monotonicity assumption:
\[\label{h}
\tag{h}\;\; \begin{array}{l}
  \text{There exists an $R_0 >0$ such that the map $r \mapsto W(a^\pm
  +r\xi)$ has a strictly}\\
  \text{positive derivative for every $r \in
  (0,R_0)$ and every $\xi \in \R^N, |\xi|=1$.}
\end{array}
\]
This is a rather weak non-degeneracy assumption, allowing for
potentials with degenerate \emph{$C^{\infty}$- flat} minima. From
now on we assume, as we can, that $r_0 <R_0$, hence
$\mathbb{B}(a^{\pm},r_0)$ are in the monotonicity region. We will
need to express $U_L$ in polar form: for any $U$ in
$[H^1_{\textrm{loc}}\weight$, we set $U^{\pm}(x):= a^{\pm} +
\rho^{\pm}(x) n^{\pm}(x)$. Then $|(U^{\pm})_x|^2 =
 ((\rho^{\pm})_x)^2 + (\rho^{\pm})^2 |n^{\pm}_x|^2$. For any $I \subseteq \R$ measurable,
we shall interpret integrals expressed in polar form as
\[
\int_I \big|U_x\big|^2 e^{cx}dx = \int_{I \bigcap \{\rho^{\pm}>0\}}
\Big\{(\rho^{\pm}_x)^2 + (\rho^{\pm})^2 |n^{\pm}_x|^2
\Big\}e^{cx}dx,
 \]
\noi since the imbedding $[H^1_{\textrm{loc}}(\R)]^N \hookrightarrow
[C^0_{\textrm{loc}}(\R)]^N$ implies $|U_x| =0$ a.e.\ on the closed
sets $U^{-1}(\{a^{\pm}\})$, even when they have positive measure.
For any $\mu < \nu$ in $\R$, we set
\[
E_c\big( U ,({\mu},{\nu})\big) :=
\int_{\mu}^{\nu}\Big\{\frac{1}{2}\big|U_x \big|^2 + W(U)
\Big\}e^{cx}dx.
\]
This is the action (\ref{action}) restricted on $[\mu, \nu]$.

\bl (cf. \cite{A-F}) \label{Interior Lemma} Assume W satisfies $(h)$
and $c>0$ is fixed. Let $a \in \{a^+, a^-\}$ and $U \in [H^1(\mu,
\nu)]^N$ with $U=a+ \rho n$, and suppose that

\noindent (i) $0< \rho(\mu)= \rho(\nu)=r \leq R_0$ ($R_0$ as in
\eqref{h}),

\noindent (ii) $r \leq \rho(x) \leq R_0$, for all $x \in (\mu,
\nu)$.

\noindent Then, there exists a $\widetilde{{U}} \in [H^1(\mu,
\nu)]^N$, $\widetilde{U}=a+ \widetilde{{\rho}} n$, such that
$U(\mu)=\widetilde{{U}}(\mu)$, $U(\nu)=\widetilde{{U}}(\nu)$ and
$\widetilde{\rho}(x) < r$, for all $x \in (\mu, \nu)$ while
\[
E_c(\widetilde{{U}},({\mu},{\nu})) \;<\; E_c({U},({\mu},{\nu})).
\]
\el \noindent In particular, locally minimizing solutions to $\eq$
on $[\mu, \nu]$ attain the maximum value $r$ of their polar radius
$\rho^{\pm} = \big| U - a^{\pm} \big|$ only at the endpoints $\{ \mu
\}, \{ \nu \}$.

\BPL \ref{Interior Lemma}. We note that the proof of Lemma 3.3 in
\cite{A-F} is based on a pointwise deformation and thus it holds
generally for functionals of the form $\int(\frac{1}{2}|U_x|^2
+W(U))d \mu(x)$ with $\mu$ a positive Radon measure. See Lemma
\ref{l1} for a similar argument. \qed

\noindent We now prove that in view of \eqref{h}, the polar radii of
$U_L$ are weak subsolutions of the operator $L(\rho):= \rho_{xx} + c
\rho_x$ in $[H^1(\mu, \nu)]^N$, for all $\mu < \nu < -L$ and $L <
\mu < \nu$.  We write $U^{\pm}(x):= a^{\pm} + \rho^{\pm}(x)
n^{\pm}(x)$ (cf. Stefanopoulos \cite{Stef}).

\begin{proposition} \textbf{(Constrained minimizers as radially
weak $H^1$ subsolutions)} \label{Constraint minimizers as radially
weak $H^1$ subsolutions} The minimizers $U_L$ of Theorem
\ref{Constrained Problem} satisfy
\[
-\big(\rho^{\pm}_L \big)_{xx}- c \big(\rho^{\pm}_L \big)_x +
\rho^{\pm}_L \big|(n^{\pm}_L)_x \big|^2+ \nabla W\big(a^{\pm} +
\rho^{\pm}_L n^{\pm}_L\big) \cdot n^{\pm}_L \leq 0,
\]
weakly in $H^1_{\textrm{loc}}\big((L, \infty) \bigcap
\{\rho^{+}_L>0\}\big)$ and $H^1_{\textrm{loc}}\big((-\infty, -L)
\bigcap \{\rho^{-}_L>0\}\big)$. In particular, if $W$ satisfies
\eqref{h}, we obtain
\[
\big(\rho^{\pm}_L \big)_{xx} + c \big(\rho^{\pm}_L \big)_x \geq 0.
\]
\end{proposition}

\BPP \ref{Constraint minimizers as radially weak $H^1$
subsolutions}. We construct local variations that do not violate the
constraint $\rho^{\pm}_L  \leq r_0$. For definiteness we consider
the case $a=a^-$, the other is similar. We take
$\phi(x):=\theta(x)n_L^-(x)$, with $\theta$ in $C^{\infty}_c
(-\infty, -L) $ and consider one-sided variations of the form
\[
U_L^{\e}(x):=U_L(x) - \e \phi(x) = a^- + \big(\rho^-_L(x)- \e
\theta(x)\big)n^-_L(x)
\]
(radially inside the cylinder) which satisfy the constraint for
small $\e \in [0, \e_{\phi}]$. Since $U_L$ is a minimizer of $E_c$,
$E_c(U_L^{\e}) \geq E_c(U_L)$, for all $\e \in [0, \e_{\phi}]$.
Consequently,
\[
\underset{\e \rightarrow 0^+}{\lim}\Big[\frac{1}{\e}\Big(
E_c(U_L^{\e}) - E_c(U_L)\Big) \Big] \geq 0.
\]
We calculate, using that $\textrm{supp}(\theta) \subseteq (-\infty,
-L)$,
\begin{align*}
E_c(U_L^{\e})\ =& \  \int_{-\infty}^{-L}\Bigg\{ \frac{1}{2}\big(
(\rho_L^-)_x - \e \theta_x\big)^2 \ +\ \frac{1}{2}\big( \rho_L^- -
\e
\theta \big)^2 \big|(n_L^-)_x \big|^2\\
&\ + \ W\big(a^- + (\rho_L^- - \e
\theta) n_L^-\big)\Bigg\}e^{cx}dx\\
& \  +\; \int^{\infty}_{- L}\left\{ \frac{1}{2}\big|(U_L)_x \big|^2
+ W(U_L) \right\}e^{cx}dx.
\end{align*}
Taking one-sided $\frac{d}{d \e}\big|_{\e = 0^+}$, we get
\[
\int_{-\infty}^{-L}\left\{-(\rho_L^-)_x(\theta_x e^{cx}) -
\Big[\rho_L^- \big|(n_L^-)_x \big|^2 + \nabla W\big(a^- + \rho_L^-
n_L^-\big)\cdot n_L^- \Big] (\theta e^{cx})\right\} dx \;\geq\; 0.
\]
We write $\theta_x e^{cx}= (\theta e^{cx})_x - c\theta e^{cx}$ and
substitute to get
\[
\int_{-\infty}^{-L}\Bigg\{\big(\rho_L^- \big)_x \big(\theta
e^{cx}\big)_x + \Big[\rho_L^- \big|(n_L^-)_x \big|^2 - \nabla
W\big(a^- + \rho_L^- n_L^-\big)\cdot n_L^- \Big] (\theta
e^{cx})\Bigg\} dx \;\leq\; 0.
\]
We are done, since the multiplication operator $M_{e^{c \text{Id}}}$
is a Fr\'echet automorphism on the dense subspace $C^{\infty}_c
\big(-\infty, -L\big)$ of $H^1_{\textrm{loc}}\big(-\infty, -L\big)$.
\qed

\noindent It is now straightforward that all $(U_L)_{L \geq 1}$
realize the constraint at most at the rims of the cylinders.

\begin{proposition} \textbf{(Contact at most at the rims of the cylinders)}
\label{Contact at most at the rims} If $W$ satisfies \eqref{h}, then

\noindent a) If $x_L^+:=\inf\big\{t \in \R : \rho^{+}_L\leq r_0
\text{ on }[t,+\infty)\big\}$, then we have $\rho^{+}_L < r_0$ on
$(x^+_L, +\infty)$.

\noindent b) If $x_L^-:=\sup\big\{t \in \R : \rho^{-}_L\leq r_0
\text{ on }(-\infty,t]\big\}$, then we have $\rho^{-}_L < r_0$ on
$(-\infty, x^-_L)$.
\end{proposition}

\BPP \ref{Contact at most at the rims}. We drop sub/superscripts $L,
\pm$ for $\rho$ and prove only a), since b) is analogous. By
definition, $x_L^+ \in (-L,L]$ and it is the time at which $U_L$
enters $\mathbb{B}(a^+,r_0) $ and remains inside it for all later
times. Minimizers $U_L$ are, by \eqref{h}, radially weak $H^1$
subsolutions: $\rho_{xx} + c \rho_x \geq 0$. Let $x_0 \in (x^+_L
,\infty)$ be such that $\rho(x_0)=r_0$. Since the point $x_0$ lies
in the interior of $[x^+_L, x_0 +1]$, by the Strong Maximum
Principle for weak $C^0$ subsolutions (\cite{G-T}), we have that
either $\rho(x_{0})<r$, or $\rho \equiv r_0$ on $[x^+_L, x_0 +1]$.
Lemma \ref{Interior Lemma} implies that $\rho$ is not identically
$r_0$, otherwise we obtain a contradiction to minimality of $U_L$.
Hence, $\rho < r_0$ on $(x^+_L,+\infty)$. \qed

\begin{proposition} \textbf{(Constrained minimizers are piecewise solutions)}
\label{Constraint minimizers as piecewise solutions} All $U_L$ are
solutions to $\eq$ in $[C^2_{\textrm{loc}}(\R \set \{x^\pm_L\})]^N
\bigcap [C^0(\R)]^N$. They are in $[C^2_{\textrm{loc}}(\R)]^N$
except possibly when $x^\pm_L = \pm L$.
\end{proposition}

\BPP \ref{Constraint minimizers as piecewise solutions}. By
Proposition \ref{Contact at most at the rims}, $|U_L(x)-a^{\pm}| <
r_0$, for all $x \in \R \set [x^-_L,x^+_L]$. Take any point $x^* \in
\R \set \{x^-_L,x^+_L\}$. By continuity, there exists an $\e_0 >0$
and a compact tubular neighborhood $\big\{\mathbb{B}(U_L(x),\e_0)\;
: \;x \in [x^* - \de, x^* + \de]\big\}$ of the graph of $U_L$ not
intersecting the boundary of the constraint cylinders, the assertion
being trivial when $x^* \in (x^-_L,x^+_L)$. This holds for $x^\pm_L$
as well, when $x^+_L <L$ and $x^-_L>-L$. We take variations of $U_L$
the $U^{\e}_L:=U_L - \e \phi$, $|\e| \leq \e_0$ small, for all
$\phi$ in $[C_c^{\infty}(x^* - \de -\e_1, x^* + \de+\e_1)]^N$, $\e_1
>0$ small, whose restriction on $(x^* - \de, x^* + \de)$ is dense in
$[H^1(x^* - \de, x^* + \de)]^N$. Using that $\phi_x e^{cx}=(\phi
e^{cx})_x - c \phi e^{cx}$, we easily get that $\eq$
 is solved weakly. Since $\nabla W \in [C^1_{\textrm{loc}}(\R^N)]^N$ and $(U_L)_x \in
 [L^2_{\textrm{loc}}(\R)]^N$, there exists $(U_L)_{xx} \in
 [L^2_{\textrm{loc}}(\R \set \{x^\pm_L\})]^N$ and therefore $U_L
 \in [C^1_{\textrm{loc}}(\R \setminus \{x^\pm_L\})]^N$ which gives that $U_L
 \in [C^2_{\textrm{loc}}(\R \setminus \{x^\pm_L\})]^N$, since $\nabla W \in
 [C^1_{\textrm{loc}}(\R^N)]^N$.   \qed

\begin{remark} (i) \textbf{(Polar form of the equation)} Write the equation $\eq$ in polar coordinates $U_L=a^{\pm} +
\rho_L^{\pm} n_L^{\pm}$ and multiply by $n_L^{\pm}$ to get that the
polar radii $\rho_L^{\pm}$ of $U_L$ satisfy the equation
\begin{equation}
 \label{Equation in polar}
 \;\;\;\; (\rho)_{xx} \;+\;c (\rho)_{x} \;=\;
\rho |n_x |^2 \;+\; \nabla W(a^{\pm} +\rho n) \cdot n.
\end{equation}
(ii) \textbf{(Energy formula)} Integrating once the equation as in
the proof of Lemma \ref{Characterization of c}, we get the formula
 \be \label{action formula}
 c \int_{\mu}^{\nu}\big| U_x \big|^2 dx = \Bigg( W(U) -
 \frac{\big| U_x \big|^2}{2}\Bigg) \Bigg|_{\mu}^{\nu},
\ee
 on any interval $[\mu, \nu]$, on which $U$ solves $\eq$
classically.
\end{remark}
\begin{proposition} \textbf{(Asymptotic behavior of constrained
minimizers)} \label{Asymptotic behavior of constraint minimizers} If
W satisfies \eqref{h}, then $U_L(x) \larrow a^{\pm}$ as $x
\rightarrow \pm \infty$. Moreover, the polar radii $\rho^{\pm}_L$ of
$U_L$ are eventually strictly monotone inside the cylinders and also
$(U_L)_x( \pm \infty) = 0$ at least up to sequences.
\end{proposition}

\BPP \ref{Asymptotic behavior of constraint minimizers}. We treat
both cases together, dropping indices $\pm , L$ of $\rho$.

\ms

\noindent  \textbf{Claim 1.} \emph{The polar radii are eventually
strictly monotone in the cylinders.}

\noi Indeed, by Lemma \ref{Interior Lemma} and the action minimality
of $U_L$, $\rho$ can not be identically constant on any subinterval
of $(-\infty,x_L^-)$ $(x_L^+ ,\infty)$. Hence, by continuity of
$\rho$ the set of critical points $A:=\{\rho_x=0\}$ is discrete.
Since $\rho$ solves $\rho_{xx}+c\rho_x \geq0$, the Maximum Principle
implies that $A$ does not contain maximum points. Moreover, $A$ can
not contain more than one minimum point; if a minimum point exists,
then at all latter points (in the unbounded direction of time)
$\rho_x$ preserves its sign on both sides of the critical point.
Hence, $\rho$ is eventually strictly monotone.

\ms

\noindent Let now $r^*$ denote the asymptotic limit of $\rho$. At $+
\infty$ it readily follows that $r^*=0$, since $e^{c
\text{Id}}W(U_L)$ is in $L^1(L,\infty)$. Indeed,
\begin{align*}
\int_L^{\infty} W^+(U_L)e^{cx}dx\;& \leq \; E_c(U_{aff})\;+\;
\int_{-\infty}^L W^-(U_L)e^{cx}dx \\
 & \leq \; E_c(U_{aff})\;+\;
\frac{W^-(a^-)e^{cL}}{c} \;< \;\infty
 \end{align*}
and $a^+$ is the only zero of $W$ inside the ball
$\mathbb{B}(a^+,r_0)$. Now we consider the limit at $-\infty$.

\noi

\noindent  \textbf{Claim 2.} \emph{For any $t\in \R$ such that
$[t,t+1]\sub (-\infty,x_L^-)$, we have
 \be \label{e14}
0 \ \leq \ \underset{|\xi|=1}{\underset{t \leq s \leq t+1}{\min}}
\Big[\nabla W\big(a^- + \rho(s) \xi\big) \cdot \xi\Big] \; \leq \;
\rho_x(t+1) e^{c}\ - \ \rho_x(t).
 \ee}

\noi Indeed, since $\eq$ is solved by $U_L$ on $(-\infty, x^-_L)$,
we integrate once the $e^{cx}$ - multiple of equation (\ref{Equation
in polar}) on $[t, t+ 1]$ to find
\begin{align*}
\int_t^{t+1}\big(\rho_x e^{cx} \big)_x dx \;& =\; \int_t^{t+1}e^{cx}
\Big(\nabla W(a^- + \rho n) \cdot n + \rho \big|n_x \big|^2\Big)
dx \\
& \geq \;e^{ct} \int_t^{t+1} \Big(\nabla W(a^- + \rho n) \cdot n +
\rho \big|n_x \big|^2\Big) dx\\
& \geq \;e^{ct} \int_t^{t+1} \nabla W(a^- + \rho n) \cdot n\
dx\\
& \geq \; e^{ct} \underset{s \in [t,t+1]}{\min} \Big[\nabla W(a^- +
\rho(s) n(s)) \cdot n(s)\Big]\\
& \geq \;  e^{ct}  \underset{|\xi|=1}{\underset{t \leq  s \leq
t+1}{\min}} \Big[\nabla W\big(a^- + \rho(s) \xi\big) \cdot \xi\Big].
\end{align*}
Utilizing assumption \eqref{h}, we obtain \eqref{e14}.

\medskip

\noi Since the limit of $\rho$ at $- \infty$ exists, there exists a
sequence $x_n \larrow -\infty$ such that $\rho_x(x_n) \larrow 0$.
Suppose first that eventually $\rho_x \geq 0$. By setting $t:=x_n
-1$ in \eqref{e14} and employing the monotonicity of $\rho$, we have
\[
0 \ \leq \ \underset{|\xi|=1}{\min} \Big[\nabla W\big(a^- + \rho(x_n
-1) \xi\big) \cdot \xi\Big] \; \leq \; \rho_x(x_n) e^{c}.
\]
By employing that $\rho(x_n -1)\larrow r^*$ and that $\rho_x(x_n)
\larrow 0$ as $n \larrow \infty$, in the limit we obtain $\nabla
W\big(a^- + r^* \xi\big) \cdot \xi=0$ for some $\xi$. Since $a^-$ is
the only critical point in $\mathbb{B}(a^-,r_0)$, it follows that
$r^*=0$. Similarly, if $\rho_x \leq 0$, we take $t:=x_n$ to get
\[
0 \ \leq \ \underset{|\xi|=1}{\min} \Big[\nabla W\big(a^- +
\rho(x_n) \xi\big) \cdot \xi\Big] \; \leq \; \big| \rho_x(x_n) \big|
\]
and again by passing to the limit as $n \larrow \infty$ it follows
that $r^*=0$.

\noi Now we consider the convergence of the derivative. By
multiplying (\ref{Equation in polar}) by $\rho$ and adding
$(\rho_x)^2$, we obtain the identity
\begin{equation}
 \label{Identity polar equation}
\big|U_x \big|^2 \;+\;\rho \nabla W(U) \cdot n \;=\; \frac{1}{2}
\Big[(\rho^2)_{xx}\;+ c\; (\rho^2)_{x} \Big].
\end{equation}
Since $\rho^2$ is also strictly increasing and has a limit at
$-\infty$, we get $(\rho^2)_x \geq 0$ and that there exists a
sequence $\xi_n \larrow -\infty$ such that $(\rho^2)_x(\xi_n)
\larrow 0$. By (\ref{Identity polar equation}), assumption \eqref{h}
and integration on $[\xi_n -1, \xi_n]$, we get
\begin{align*}
0 \;\leq \; \int_{\xi_n -1}^{\xi_n}\big|U_x \big|^2 dx\; & \leq \;
\frac{1}{2} \Big[(\rho^2)_{x}(\xi_n)-(\rho^2)_{x}(\xi_n -1)\Big]\;+
\frac{c}{2}\Big[\rho^2(\xi_n)-\rho^2(\xi_n-1) \Big]\\
& \leq \; \frac{1}{2} \Big[(\rho^2)_{x}(\xi_n)\;+\;c\; \rho^2(\xi_n)
\Big] \;\larrow 0,
\end{align*}
as $n \larrow \infty$. The proof is complete. \qed

\ms

\noindent We conclude this section by proving that $(U_L)_x \in
[L^2(\R)]^N$, but not $L$-uniformly. In addition, $U_L$ satisfies
the first formula of Lemma \ref{Characterization of c}
approximately, up to some additional terms which \emph{relate $c$
with the jump of $(U_L)_x$ at the rims}.

\begin{proposition} \textbf{(Approximate relation for $c$)}
\label{Approximate relation for c} The 1-sided derivatives
$(U_L)_x(\pm L^{\pm})$ of $U_L$ exist, and
\begin{align*}
c \int _{\R} \big|(U_L)_x\big|^2 dx \;=\; W^-(a^-) \; &+
\;\frac{1}{2}\Big(\big|(U_L)_x(-L^{+})\big|^2 -\big|(U_L)_x(
-L^{-})\big|^2\Big) \\
& +\; \frac{1}{2} \Big(\big|(U_L)_x(+L^{+})\big|^2
-\big|(U_L)_x(+L^{-})\big|^2\Big).
\end{align*}
In particular, $(U_L)_x \in [L^2(\R)]^N$.
\end{proposition}

\BPP \ref{Approximate relation for c}. Proposition \ref{Constraint
minimizers as piecewise solutions} assures that we can apply formula
(\ref{action formula}) on $(-\infty, -L- \e)$, $(-L+\e, L- \de)$ and
$(L+\de, \infty)$ for $\e, \de>0$ small utilizing by \ref{Asymptotic
behavior of constraint minimizers} the asymptotic behavior of
$U_L$'s and the continuity of $W$. We obtain three relations on
these intervals. Utilizing H\"older's inequality, we easily find
\begin{align*}
\big|(U_L)_x(-L-\e)\big| & \leq \sqrt{2} \Big( W(U_L(-L-\e)) \;+\;
W^-(a^-)\Big)^{\frac{1}{2}} ,\\
 \big|(U_L)_x(-L+\e)\big| &\leq
\sqrt{2} \Bigg( ce^{c(L- \e)} \int_{-L+\e}^{L-\de}
\big|(U_L)_x\big|^2e^{cx}dx \;-\; W(U_L(+L-\de)) \\
& \ \ +\; W(U_L(-L+\e)) \;+\; \frac{1}{2}
\big|(U_L)_x(+L-\de)\big|^2\Bigg)^{\frac{1}{2}} ,\\
\big|(U_L)_x(+L-\de)\big| & \leq \sqrt{2} \Bigg( W(U_L(L-\de))-
W(U_L(-L+\e)) + \frac{1}{2}
\big|(U_L)_x(-L+\e)\big|^2\Bigg)^{\frac{1}{2}} ,\\
\big|(U_L)_x(+L+\de)\big| & \leq \sqrt{2} \Bigg( ce^{-c(L+
\de)}\int^{\infty}_{L+\de} \big|(U_L)_x\big|^2e^{cx}dx \;+\;
W(U_L(L+\de)) \Bigg)^{\frac{1}{2}}.
\end{align*}
Letting $\e \larrow 0^+$ and $\de \larrow 0^+$ separately, we obtain
that the moduli of the one-sided limits exist, but may differ.
Adding these relations and letting $\e, \;\de \larrow 0^+$ we obtain
the formula for $c$. \qed

\section{The Local Replacement Lemmas.}

We recall some basics from Differential Geometry. The canonical
coordinates $(p,d)$ on $\R^N$ with respect to a $C^2$ convex set
$\mC \sub \R^N$ are defined by
 \be \label{e1}
 u\ =: \ p \ + \ d n
 \ee
where $p$ is the projection on the convex set $\mC$, $0\in \mC$, $d$
the signed distance from $\p \mC$ and $n$ the outward unit normal of
$\p \mC$. The latter is parameterized by the $C^2$ local coordinates
\[
\R^{N-1} \ni s\ =\ (s_1,...,s_{N-1})\ \mapsto \ p(s_1,...,s_{N-1})
\in \p \mC.
\]
We may assume that the set of vectors
 \be \label{e2}
\frac{\p p}{\p s_i}\ = \ \vec{t_i}\ , \ \ \ \ i=1,...,N-1,
 \ee
is an orthonormal frame in the tangent space at $p$, coinciding with
the principal curvature directions (\cite{DC}, p.\ 144, p.\ 216).
Thus,
 \be \label{e3}
\frac{\p n}{\p s_i}\ = \ \kappa_i \vec{t_i}, \ \ \
\kappa_i=\kappa_i(s) \text{ the i-th principal curvature of }\p \mC.
 \ee
The coordinate system $(p,d)$ is defined for $-d_0\leq d$, provided
that $d_0 \kappa_i \leq 1$, $i=1,...,N-1$ (\cite{G-T}). The
orientation is such that $\kappa_i\geq0$ when $\mC$ is convex. We
write
 \be \label{e4}
U(x)\ = \ p(x) \ +\ d(x) n(x),
 \ee
meaning $p(x)=p(s(x))$, $n(x)=n(s(x))$. By differentiating
\eqref{e4},
 \begin{align*}
\dot{U}(x)\ &= \ \dot{p}(x)\ + \ \dot{d}(x)n(x)\ + \ d(x)\dot{n}(x)\\
            &= \ \vec{t_i}\dot{s}_i \ + \ \dot{d}n \ + \ d \kappa_i
            \vec{t_i} \dot{s}_i.
 \end{align*}
Hence,
 \be \label{e5}
|\dot{U}(x)|^2\ = \ \sum_{i=1}^{N-1} \dot{s}_i^2(1+\kappa_i d(x))^2
\ +\ (\dot{d}(x))^2.
 \ee
Let now $\mC' \sub \R^N$ be a convex set and assume that
 \be \label{e6}
W_u \cdot n \ \geq \ \frac{c_0}{2} >0\ \ \ \text{on }\p \mC',
 \ee
where $W \in C^1(\R^N)$ and $(p,d)$ the canonical coordinates
associated to $\p \mC'$. By the $C^1$ smoothness of $W$ and
\eqref{e6}, there is a $\bar{d}>0$ such that
  \be \label{e7}
d \ \lmapsto \ W(p+dn) \ \ \ \text{is increasing for } -\bar{d} \leq
d \leq \bar{d}.
  \ee

\bl \label{l1} Let $x_1<x_2$ in $\R$ and $U \in [H^1(x_1,x_2)]^N$ be
such that

\noi (i) $d(x_1)=d(x_2)=0$,

\noi (ii) $0\leq d(x) \leq \bar{d}$, for $x\in(x_1,x_2)$.

\noi If \eqref{e6} and \eqref{e7} are satisfied, then there exists
$\tilde{U} \in [H^1(x_1,x_2)]^N$ with the following properties:
\[
\tilde{U}(x_1)\ = \ {U}(x_1),\ \ \ \tilde{U}(x_2)\ = \ {U}(x_2),
\]
\[
-\bar{d} \ \leq \ \tilde{d}(x)\ <\ 0, \text{ for $x\in(x_1,x_2)$},
\]
\[
E_\mu(\tilde{U},(x_1,x_2))\ < \ E_\mu(U,(x_1,x_2)),
\]
where $\tilde{U}(x)=\tilde{p}(x)+\tilde{d}(x)n(x)$ and
\[
E_\mu(U,(x_1,x_2))\ :=\
\int_{x_1}^{x_2}\left(\frac{1}{2}|\dot{U}(x)|^2+W(U(x))\right)d\mu(x)
\]
where $\mu$ is a positive Radon measure on $\R$. \el

\BPL \ref{l1} (cf. Lemma 3.3 in \cite{A-F}). Let $\phi:[0,1]\larrow
\R$ be a smooth function such that $\phi(0)=\phi(1)=0$,
$\phi(\sigma)>0$ for $\sigma \in (0,1)$. For small $\e \geq 0$
define
\[
\widetilde{U}^\e (x)\ :=\ p(x) - \e
\phi\left(\frac{x-x_1}{x_2-x_1}\right)n(x), \ \ \ x \in [x_1,x_2],
\]
where $U(x)=p(x)+d(x)n(x)$. By \eqref{e5}, we have
\[
|\dot{U}(x)|^2\ = \ \sum_{i=1}^{N-1} \dot{s}_i^2(x) \ + \ d^2
\sum_{i=1}^{N-1} \ \kappa^2_i \dot{s}^2_i(x) \ +\ 2d
\sum_{i=1}^{N-1} \kappa_i \dot{s}_i^2(x) \ + \ \dot{d}^2(x).
\]
We note that
\[
|\dot{\tilde{U}}^\e|^2\ = \ \sum_{i=1}^{N-1} \dot{s}_i^2 \ + \
\e^2\phi^2 \sum_{i=1}^{N-1} \kappa^2_i \dot{s}^2_i \ - \ 2\e \phi
\sum_{i=1}^{N-1}  \kappa_i \dot{s}^2_i \ + \ \e^2
\frac{\phi'^2}{(x_2 -x_1)^2}.
\]
Thus, we have that
 \begin{align} \label{e9}
E_\mu(\tilde{U}^\e,(x_1,x_2))\ =&  \ E_\mu(\tilde{U}^0,(x_1,x_2))\nonumber\\
 &-\ \e \int_{x_1}^{x_2}\phi \sum_{i=1}^{N-1}  \kappa_i \dot{s}^2_i
 d\mu \ + \ \frac{\e^2}{2} \int_{x_1}^{x_2}\phi^2 \sum_{i=1}^{N-1}  \kappa_i^2 \dot{s}^2_i
 d\mu\\
 &- \ \int_{x_1}^{x_2}\Big(W(p)\ - \
 W(p-\e\phi n)\Big)d\mu \nonumber\\
 &+\ \frac{\e^2}{(x_2 - x_1)^2}\int_{x_1}^{x_2}\phi'^2d\mu \nonumber.
 \end{align}
By \eqref{e7}, (ii) above and convexity of $\mC'$ we have
 \be \label{e8}
E_\mu(\tilde{U}^0,(x_1,x_2))\ \leq  \ E_\mu(U,(x_1,x_2)).
 \ee
On the other hand, \eqref{e7} also implies
 \begin{align}
- \ \int_{x_1}^{x_2} \Big(W(p)\ -& \
 W(p-\e\phi n)\Big)d\mu \
 +\ \frac{\e^2}{2(x_2 - x_1)^2} \int_{x_1}^{x_2}\phi'^2 d\mu &
  \nonumber\\
 =&\
 - \int_{x_1}^{x_2} \left(\int_0^1 \frac{d}{d\tau}(W(p-\e \tau \phi n ))
 d\tau\right) d\mu\nonumber\\
 &+\  \frac{\e^2}{2(x_2 - x_1)^2} \int_{x_1}^{x_2} \phi'^2 d\mu
\\
  =& \ - \e \int_{x_1}^{x_2} \left( \int_0^1 W_u (p -\e \tau \phi n) \cdot
  \phi n
 \right) d\tau d\mu \nonumber\\
 &+\  \frac{\e^2}{2(x_2 - x_1)^2} \int_{x_1}^{x_2}\phi'^2 d\mu
  \nonumber\\
 \overset{\eqref{e6}}{<}& \ - C \e \ + \ \frac{\e^2}{2(x_2 - x_1)^2}
\int_{x_1}^{x_2} \phi'^2 d\mu \ < \ 0,  \nonumber
 \end{align}
for some $C>0$ and small $\e>0$. Finally, we observe that by the
convexity of $\mC'$,
\[
-\ \e \int_{x_1}^{x_2}\phi \sum_{i=1}^{N-1}  \kappa_i \dot{s}^2_i
 d\mu \ + \ \frac{\e^2}{2} \int_{x_1}^{x_2}\phi \sum_{i=1}^{N-1}  \kappa_i^2 \dot{s}^2_i
 d\mu  \leq \ 0,
\]
for small $\e>0$. From these inequalities and \eqref{e9}, the lemma
follows with $\tilde{U}:=\tilde{U}^\e$, $0<\e<<1$. \qed

\ms

\noi \textbf{Hypotheses}

\ms

\noi (H1) $W:\R^N \larrow \R$, $C^2$, with two minima
$W(a^-)<W(a^+)=0$.

\ms

\noi (H2) $\{u|W(u)\leq0\} =: \mC^-_0\cup \{a^+\}$, $\mC^-_0$
compact, convex.

\ms

\noi (H3) (i) $W_u \cdot n \geq c_0 >0$ on $\p \mC^-_0 =:
\{W=0\}_{(-)}$, $n$ the outward unit normal on $\p \mC_0^-$.

\ \ \ (ii) $W_{uu}\geq c_0 I$ on $\{W=0\}_{(-)}$.

\begin{remark} a) By $C^2$ smoothness of $W$, there exists a $b>0$ such
that
 \be \label{e10}
W_{uu}\ \leq \ bI, \ \ \ \text{on } \{u|W(u)\leq0\}.
 \ee
\noi b) (H3) implies that the set $\{u|W(u)=\beta\}$ for
$0<\beta<<1$ is made up of two components, which we denote by
\[
\{W=\beta\}_{(-)}\ \ \ \text{and}\ \ \ \{W=\beta\}_{(+)},
\]
with $\{W=\beta\}_{(-)}$ convex and enclosing $a^-$. On the other
hand, for $\beta<0$ ($|\beta|<<1$), $\{u | W(u)=\beta\}$ is made up
of one component which is convex. So more precisely there is an
$\al_0>0$ such that $\{W=\beta\}_{(-)}$ is convex,
$\al_0\leq\beta\leq\al_0$. By the smoothness of $W$,
 \be \label{e11}
W_u \cdot n \ \geq \ \frac{c_0}{2} \ \ \ \text{on
$\{W=\beta\}_{(-)}$, $\al_0\leq\beta\leq\al_0$}.
 \ee
Note that the sets $\{W=\beta\}_{(-)}$ are nested for
$\al_0\leq\beta\leq\al_0$.
\end{remark}

\noi Now we take $\al \in (0,\al_0)$ and furthermore restrict it as
follows:
 \be \label{e12}
0\ < \ \al\ < \frac{c_0}{4}\la\ =:\ \bar{a}_0,
 \ee
where $\la$ is a fixed number satisfying the conditions
\[
0\ \leq \ \la \ \leq \ \frac{c_0}{2b}\ , \ \ 0\ <\ \la \ \leq \ d_0,
\ \ \la \ < \ \frac{1}{\max\{\kappa_1,...,\kappa_{N-1}\}},
\]
with $b$ as in \eqref{e10} above,
\[
d_0\ = \ \dist\Big(\{W=\al_0\}_{(-)},\{W=-\al_0\}_{(-)}\Big),
\]
and $\kappa_1,...,\kappa_{N-1}$ the principal curvatures of
$\{W=\beta\}_{(-)}$ (all positive by convexity). We note that
 \be \label{e13}
W(p\ - \ \la n(p))\ < \ 0, \ \ \ \text{for } p \in \{W=\al\}_{(-)}.
 \ee
Indeed (dropping $p$ in $n(p)$),
\begin{align*}
W(p)\ - \ W(p-\la n) \ &=\ -\int_0^\la \frac{d}{dt}[W(p-tn)]dt\\
 &=\ \int_0^\la (W_u(p-tn)-W_u(p)+W_u(p))\cdot n \ dt\\
 &=\ \int_0^\la W_u(p)\cdot n dt \ - \  \int_0^\la \int_t^0
   \frac{d}{ds}(W_u(p-sn)ds)\cdot n \ dt\\
 &=\ \int_0^\la W_u(p)\cdot n dt \ - \  \int_0^\la \int^t_0
   W_{uu}(p-sn)n\cdot n \ ds dt\\
 &\geq \ \frac{c_0}{2}\la \ - \ \frac{b}{2}\la^2 \ \ \ \ \ \ \  (\eqref{e10},\eqref{e11})\\
 & \geq \ \frac{c_0}{4}\la  \ \ \ \ \ \ \ \ \ \ \ \ \ \ \ \ \ \ \  \left(\la\leq\frac{c_0}{2b}\right).
\end{align*}
Therefore, we have
\[
W(p)-\frac{c_0}{4}\la \ \geq \ W(p-\la n)
\]
and so (by \eqref{e12})
\[
0\ > \ \al \ - \ \frac{c_0}{4}\la \ \geq \ W(p-\la n).
\]

\bl \label{l2} Let $\mC$ denote the component of $\{u|W(u)\geq
\al\}$ with $\p \mC =
 \{W=\al\}_{(-)}$. Let $(p,d)$ be the canonical coordinates with
 respect to $\mC$. Assume that $\al$ is as in \eqref{e12}, and
 assume that (H1), (H2), (H3) hold. Let also $x_1<x_2 \in \R$ and
 $U \in [H^1(x_1,x_2)]^N$ be such that

 \noi (i) $d(x_1)=d(x_2)=0$,

 \noi (ii) $d(x_0) \geq 0$, for some $x_0 \in (x_1,x_2)$.

\ms

\noi Then, there is a $\tilde{U} \in [H^1(x_1,x_2)]^N$ with the
properties
\[
\tilde{U}(x_1)\ = \ U(x_1)\ ,\ \  \tilde{U}(x_2)\ = \ U(x_2),
\]
where $-d_0 \leq \tilde{d}(x) <0$, for $x \in (x_1,x_2)$, and
\[
E_\mu(\tilde{U},(x_1,x_2))\ <\ E_\mu(U,(x_1,x_2)),
\]
where $\tilde{U}(x)=\tilde{p}(x)+\tilde{d}(x)n(x)$.
 \el

\BPL \ref{l2} (cf. Lemma 3.4 in \cite{A-F}). Let
\[
\rho_M\ := \ \max_{x\in[x_1,x_2]}d(x).
\]
We can assume that $d(x_0)=\rho_M$. We first analyze the case
$d(x_0)=\rho_M=0$. In this case we can assume that $d(x)<0$ for some
$x \in (x_1,x_0)$ ($x \in (x_0,x_2)$), since otherwise, by Lemma
\ref{l1} we can replace $U$ with a function that satisfies this
condition and has less action. From this and the continuity of $U$
it follows the existence of $\hat{x}_1 \in (x_1,x_0)$, $\hat{x}_2
\in (x_0,x_2)$, $-\frac{d_0}{2}<\hat{d}<0$, such that
$d(\hat{x}_1)=d(\hat{x}_2)=\hat{d}$ and $\hat{d}<d(x)<0$, for $x \in
(\hat{x}_1,\hat{x}_2)$. We now consider the parallel hypersurface to
$\p \mC$, parameterized by $p+\hat{d}n(p)$, $p \in \p \mC$. This is
convex, and denote it by $\p \mC'$. It can be deduced by \eqref{e11}
that condition \eqref{e6} holds on $\p \mC'$. Then we can apply
lemma \eqref{l1} on $\p \mC'$ and obtain a local replacement between
$\hat{x}_1$ and $\hat{x}_2$ and conclude that the claim of the lemma
is true if $\rho_M =0$. Therefore we can assume $\rho_M>0$. If
$0<\rho_M\leq d_0$, again we can conclude by Lemma \ref{l1} applied
to the connected component $I_0$ of the set $\{x \in
(x_1,x_2)|d(x)>0\}$ that contains $x_0$. It remains to analyze the
case $\rho_M >d_0$. We can identify $(x_1,x_2)$ with $I_0$. Let
$h:[0,d_0] \larrow [-\la,0]$, $h(\sigma)=-\la\frac{\sigma}{d_0}$,
then $h(0)=0$, $h(d_0)=-\la$. We define the deformation
\[
\tilde{U}(x)\ := \ \left\{
\begin{array}{l}
p(x) \ +\ h(d(x))n(x),\ \ \text{for } x \in [x_1,x_2], \ d(x)<d_0 \ms\\
p(x) \ - \ \la n(x), \ \ \ \ \ \ \ \ \ \ \text{for } x \in
[x_1,x_2], \ d(x)\geq d_0,
\end{array}
\right.
\]
$\tilde{U}(x_1)=U(x_1)$, $\tilde{U}(x_2)=U(x_2)$. For the kinetic
energy we have the estimates
\begin{align*}
|\dot{U}(x)|^2\ &= \ \sum_{i=1}^{N-1} \dot{s}_i^2(1+\kappa_i d(x))^2
\ + \ \dot{d}^2(x)\\
&\geq \ \sum_{i=1}^{N-1} \dot{s}_i^2(1+\kappa_i h)^2
\ + \ (h'(d))^2\dot{d}^2(x)\\
&= \ |\dot{\tilde{U}}(x)|^2,
\end{align*}
when $d(x)<d_0$, while for $d(x)\geq d_0$ we have
\begin{align*}
|\dot{U}(x)|^2\ &= \ \sum_{i=1}^{N-1} \dot{s}_i^2(1+\kappa_i d(x))^2
\ + \ \dot{d}^2(x)\\
&\geq \ \sum_{i=1}^{N-1} \dot{s}_i^2(1+\kappa_i d(x))^2\\
&> \ \sum_{i=1}^{N-1} \dot{s}_i^2(1-\la \kappa_i )^2\\
&= \ |\dot{\tilde{U}}(x)|^2.
\end{align*}
Hence,
\[
\int_{x_1}^{x_2}|\dot{\tilde{U}}(x)|^2d\mu(x)\ < \
\int_{x_1}^{x_2}|\dot{{U}}(x)|^2d\mu(x).
\]
For the potential energy we have the estimates
\begin{align*}
W(\tilde{U}(x))\ &= W(p(x)\ + \ h(d(x))n(x))\\
&= \ W\Big(p(x)\ - \ \frac{\la d(x)}{d_0}n(x)\Big)\\
&\leq \ W(p(x)\ + \ d(x)n(x)) \ \ \ (\text{by} \eqref{e11})\\
& = \ W(U(x)),
\end{align*}
when $d(x)<d_0$, while for $d(x)\geq d_0$ we have by \eqref{e13} and
(H2) that
\[
W(\tilde{U}(x))\ \leq \ 0\ \leq \ W({U}(x)).
\]
Putting it all together, we have
\[
\int_{x_1}^{x_2}W({\tilde{U}}(x)) d\mu(x)\ < \ \int_{x_1}^{x_2}
W({{U}}(x)) d\mu(x).
\]
The argument so far establishes that
\[
E_\mu(\tilde{U},(x_1,x_2))\ <\ E_\mu(U,(x_1,x_2)).
\]
The proof of Lemma \ref{l2} is complete. \qed

\section{Action Properties of Minimizers}

We now show that $E_c(U_L)$ is a function of the jumps at the rims
$\big| (U_L)_x(\pm L^+) \big|^2 - \big| (U_L)_x(\pm L^-) \big|^2$,
while $E_c(U_L) = 0$ for minimizers in $[C^2(\R)]^N$ which solve
$\eq$ on $\R$. To prove this, we derive an \emph{equipartition}
relation \emph{at $+\infty$} (see \cite{A-Be-C}, \cite{A-F} and our
result Lemma \ref{Equipartition Limit of $E_c$ at infinity}). We
first need a formula for the action of solutions:

\bl \textbf{(1st integral)} \label{1st energy integral} Every
solution to $\eq$ in $[C^2(\mu,\nu)]^N$ satisfies:
\[
E_c(U,({\mu},{\nu})) = \int_{\mu}^{\nu}\Big\{\frac{1}{2}\big|U_x
\big|^2 + W(U) \Big\}e^{cx}dx \;=\; \Bigg\{\frac{e^{cx}}{c}\Bigg(
W(U)- \frac{\big|U_x \big|^2}{2}\Bigg) \Bigg\} \Bigg|_{\mu}^{\nu}.
\]
 \el

\BPL \ref{1st energy integral}. The equation $\eq$ implies $-U_{xx}
\cdot U_x + \nabla W(U) \cdot U_x = c \big|U_x \big|^2$, hence we
obtain
\[
\left( \frac{1}{2}\big|U_x \big|^2 - W(U) \right)_x \; =\; - c
\big|U_x \big|^2.
\]
Integrating by parts the $e^{cx}$ - multiple of this equation, we
get
\begin{align*}
\left\{ \frac{e^{cx}}{2}\big| U_x \big|^2 \right\}\Big|_{\mu}^{\nu}
\;& -\; \frac{c}{2} \int_{\mu}^{\nu}\big| U_x \big|^2 e^{cx} dx \; -
\; \Big( e^{cx} W(U) \Big)\Big|_{\mu}^{\nu} \;+ \; c
\int_{\mu}^{\nu} W(U) e^{cx}dx\\
& =\; - c \int_{\mu}^{\nu}\big| U_x \big|^2 e^{cx} dx.
\end{align*}
which leads to the desired formula. \qed

\bl \textbf{(The action in terms of the jumps)} \label{The action in
terms of derivatives' jump discontinuities} The minimizers $U_L$
satisfy
\begin{align*}
E_c(U_L) \; = \;  &\ \underset{\omega \rightarrow \infty}{\lim}\
\frac{e^{c\omega}}{c}\Bigg( W(U_L(\omega))-
\frac{\big|(U_L)_x (\omega)\big|^2}{2}\Bigg) \\
& \left.\begin{array}{l}+ \  \dfrac{e^{+cL}}{2c}
\left(\big|(U_L)_x(+L^{+}) \big|^2
  -\big|(U_L)_x( +L^{-})\big|^2 \right) \ms \ms\\
 + \ \dfrac{e^{-cL}}{2c}\left(\big|(U_L)_x(-L^{+})\big|^2
-\big|(U_L)_x( -L^{-}) \big|^2 \right).
\end{array}\right\} =:\ e_c(U_L)
\end{align*}
 \el
\noi The sum $e_c(U_L)$ comprises ``error terms'' which vanish if
$U_L \in [C^2_{\text{loc}}(\R)]^N$.

\BPL \ref{The action in terms of derivatives' jump discontinuities}.
First note that $E_c(U_L) = {\lim}_{\omega \rightarrow \infty}
E_c(U_L,({-\infty},{\omega}))$. Apply Lemma \ref{1st energy
integral} to $U_L$ which is a piecewise solution on $(-\infty, -L)$,
$(-L,L)$, $(L,\omega)$ and add the three relations, utilizing the
continuity of $W(U_L)$ at $\pm L$. Finally, let $\omega \larrow
\infty$.   \qed

\ms

\noindent Solutions to $U_{xx}=\nabla W(U)$ in the well-studied case
of $c=0$ satisfy an equipartition property: $2 W(U)= \big|U_x
\big|^2$. Our dissipation term $-c \big| U_x \big|^2$ forces a
similar behavior but \emph{at $+ \infty$}.

\bl \textbf{(Equipartition limit of the energy at $+ \infty$)}
\label{Equipartition Limit of $E_c$ at infinity} The minimizers
$U_L$ satisfy
\[
\underset{\omega \rightarrow
\infty}{\lim}\Bigg[\frac{e^{c\omega}}{c}\Bigg( W(U_L(\omega))-
\frac{\big|(U_L)_x (\omega)\big|^2}{2}\Bigg) \Bigg] = 0.
\]
 \el

\BPL \ref{Equipartition Limit of $E_c$ at infinity}. By the formula
(\ref{action formula}) for $\mu = \omega$, $\nu = \infty$ and
Proposition \ref{Asymptotic behavior of constraint minimizers}, we
have
\[
0 \ \leq \ c \int_{\omega}^{\infty}\big| (U_L)_x \big|^2 dx \;=\;
\frac{\big| (U_L)_x (\omega) \big|^2}{2}\;-\; W(U_L(\omega)).
\]
This gives
\begin{align*}
0\;\leq \;\frac{e^{c \omega}}{c}\Bigg( \frac{\big| (U_L)_x (\omega)
\big|^2}{2}\;-\; W(U_L(\omega)) \Bigg) \;& = \; e^{c \omega}
\int_{\omega}^{\infty}\big| (U_L)_x \big|^2 dx \\
& \leq \; \int_{\omega}^{\infty}\big| (U_L)_x \big|^2 e^{cx} dx.
\end{align*}
By Proposition \ref{Approximate relation for c}, we have $(U_L)_x
\in [L^2\weight$. Hence, letting $\omega \larrow \infty$ we are
done. \qed

\begin{corollary} \textbf{(The action measures the jump
discontinuities)} \label{The action equals derivatives' jump
discontinuities} We have that $E_c(U_L) = e_c(U_L)$, with $e_c(U_L)$
as in Lemma \ref{The action in terms of derivatives' jump
discontinuities}. In particular, $E_c(U_L) = 0$ if $U_L \in
[C^2_{\text{loc}}(\R)]^N$.
\end{corollary}

\section{Implications of the Local Replacement Lemmas. Determination of the Speed.}

\noi We first introduce our main hypothesis on the potential (cf.
(H1)-(H3) in Sec.\ 4):
\[
\label{h*} \tag{$h^*$} \begin{array}{l}\text{$W$ is in
$C^2_{\textrm{loc}}(\R^N)$, $a^{\pm}$ are minima,
$W(a^-)<0=W(a^+)$ and $\min_{\R^N}\{ W\}$}\\
 \text{$ =W(a^-)$. Moreover:}
\end{array}
\]
\begin{enumerate}
\item There is an $\al_0 >0$ such that for all $\al \in (0,\al_0]$,
  we have $W^{-1}(\{\al\})= \p \Ca^- \cup \p \Ca^+$, $\{u \in \R^N | W\leq \al\}
    = \Ca^- \cup \Ca^+$, where $\Ca^-, \Ca^+$ are disjoint compact, convex sets
     with $C^2$ boundaries, containing $a^{\pm}$
  respectively. Moreover, $W_u \cdot n \geq c_0>0$
  on $\p \mC_0^-$ and $W_{uu} \geq c_0 I$
  on $\p \mC_0^-$, $n$ the outward unit normal of $\p \mC_0^-$.

\item \label{it2} The map $r \mapsto W(a^- +r\xi)$ has a strictly positive
 derivative as long as $a^- +r\xi \in \Ca^-$, $|\xi|=1$, $r>0$.
\end{enumerate}

\ms

\noi Assumption \eqref{h*} implies ${{\liminf}}_{|u| \rightarrow
\infty} \big\{W(u)\big\} \geq \al_0$, thus $W$ satisfies
\[
W^{-1}\big([W(a^-), 0]\big)\ \subset \subset \ \R^N,
\]
which was assumed in Theorem \ref{Constrained Problem}.

\begin{definition} \label{d1} For $\al \in (0, \bar{\al}_0]$ and $L \geq 1$, we set
 \begin{align*}
\la^{-}_L \; & :=\; \sup \big\{ x \in \R \;\; : \;\;|U_L(x) - a^{-}|
=
r_0 \big\},\\
 \la^{+}_L \; & :=\; \inf \big\{ x \in \R \;\; : \;\;|U_L(x) - a^{+}| =
r_0 \big\},\\
\la^{\al-}_L \; & :=\; \sup \big\{ x \in \R \;\; : \;\; U_L(x) \in
\partial (\Ca^-) \big\}.
\end{align*}
\end{definition}

\noi We will show that $U_L$ intersects exactly once any of the sets
$\p \mathbb{B}(a^-,r_0)$, $\p \mathbb{B}(a^+,r_0)$, $\p \mC^-_\al$.
Decreasing $\al>0$ if necessary, we may assume $\mC_a^+ \sub
\mathbb{B}(a^+,r_0)$ and that $\mathbb{B}(a^+,r_0)$ is disjoint from
$\mC_a^-$.

\begin{proposition} \textbf{(Global a priori control on action minimizers)}
\label{Global control on minimizers} Assume $W$ satisfies \eqref{h}
and \eqref{h*}, $\al$ is as in Definition \ref{d1} and let $(U_L)_{L
\geq 1}$ be the family of minimizers of Theorem \ref{Constrained
Problem}. For all $L \geq 1$, we have

\noindent (I) $U_L$ exits $\Ca^-$ precisely once at
$x=\la^{\al-}_L$, that is
\[
x \in (-\infty, \la^{\al-}_L]\;\;\;\Longrightarrow\;\;\;W(U_L(x))
\leq \al.
\]
\noindent (II) The image $U_L(\R)$ restricted to $\R^N \set (\Ca^-
\bigcup \mathbb{B}(a^+,r_0))$ has only one connected component and
\[
\text{$W(U_L(x)) \geq \al$ for $x \in [\la^{\al-}_L,\la^{ +}_L]$}.
\]
\noindent (III) The image $U_L(\R)$ restricted on $\Ca^- \bigcup
\mathbb{B}(a^+,r_0)$ has precisely two connected components and
\[
\text{$W(U_L(x)) \leq \al$ for some $x \in \R$ then either $\;x \in
(-\infty, \la^{\al}_L]$, or $\;x \in [\la^{+}_L, +\infty)$}.
\]

\noindent (IV) The numbers $\la^{\pm}_L$ are well defined as the
unique times at which $U_L$ crosses the spheres $\partial
(\mathbb{B}(a^{\pm},r_0))$.

\noindent (V) The polar radii $\rho^{\pm}_L=\big|U_L - a^{\pm}
\big|$ are strictly monotone on $[\la^{+}_L, +\infty)$, $(-\infty,
\la^{\al-}_L]$ respectively.
\end{proposition}

\BPP \ref{Global control on minimizers}. 1. We first  settle
$\la_L^-$. We note that Lemma 3.4 of \cite{A-F} applies because the
local replacements in its proof are pointwise, and because
$W(a^-)\leq W(a^+)$. Thus, $\la_L^-$ is unique and half of (IV) is
established.

\noi 2. Next we settle $\la_L^{\al-}$. By applying Lemma \ref{l2},
we obtain the existence of a unique intersection of $U_L$ with $\p
\mC_a^-$, and so (I) is established.

\noi 3. We handle $\la_L^+$ as follows. Assume by contradiction that
$U_L$ intersects $\p \mathbb{B}(a^+,r_0)$ more than once. Then,
there are $x_1<x_2$ such that $U_L(x_i) \in \p \mathbb{B}(a^+,r_0)$,
$i=1,2$ and  $U_L(x_i) \not\in \mathbb{B}(a^+,r_0)$, $x_1<x<x_2$.
Since by step 2.\ above, $U_L$ can not intersect $\p \mC_\al^-$ for
those $x$'s, it follows that Lemma 3.4 in \cite{A-F} applies and
leads to a local replacement with less action and thus to a
contradiction. Thus, by step 1. above, (IV) has been established.

\noi 4. The previous arguments show that $U_L(x)$ can not exit
$\Ca^-$ before $x=\la^{\al-}_L$ and can not enter
$\mathbb{B}(a^+,r_0)$ before $x=\la^{+}_L$. Thus we have control on
the intervals for which $U_L$ is in the monotonicity regions, which
implies the $L^{\infty}$ bounds
\[
\big\|\rho^{-}_L\big\|_{L^{\infty}(-\infty,\la^{\al-}_L)}\;\leq\;
\underset{u \in \Ca^-}{\max}\big|u-a^-\big|\;\;,\;\;\;
\big\|\rho^{+}_L\big\|_{L^{\infty}(\la^{+}_L,\infty)} \;\leq\; r_0.
\]
It follows that Lemma \ref{Interior Lemma} can be applied to the
minimizers inside $\mathbb{B}(a^{\pm},r)$ with $r>r_0$ showing that
they can not be identically constant on any subinterval. By
\eqref{h*}, $\rho_L^{\pm}$ satisfy $(\rho_L^{\pm})_{xx} +c
(\rho_L^{\pm})_x \geq 0$. By the Strong Maximum Principle, both
$\rho_L^{\pm}$ can not have local maxima, thus they are strictly
monotone. By Proposition \ref{Asymptotic behavior of constraint
minimizers} it follows that the same is true for $r<r_0$, thus (V)
has been established. \qed

\begin{remark} \label{Remarks on lambda - times} We have the
ordering $-L\leq\la_L^{\al-}\leq\la_L^+$. We will prove existence by
showing that for some $L<\infty$ large, the constraint is not
realized: $-L < \la_L^{-}$ and $\la_L^{+} < L$ strictly. We define
\[
\la_L^{0-}\ :=\ \sup \{x\in \R \ : \ U_L(x) \in \p \mC_0^-\}.
\]
Note that $-L\leq\la_L^{0-}\leq\la_L^{\al-}\leq\la_L^+$.
\end{remark}

\noindent In the sequel we will need the following estimate.

\bl \label{Basic action estimate} If $dist(\Ca^-, \B(a^+,r_0)) =:
d_{\al}$, then for all $\al \in [0, \bar{\al}_{0}]$ and $L \geq 1$,
we have
\[
E_c(U_L) \; \geq \; - \frac{W^-(a^-)}{c}e^{c \la_L^{0-}} \;+\;
\frac{\al}{c} \Big[e^{c \la_L^{+}} - e^{c \la_L^{\al-}} \Big] \;+\;
\frac{c\; d_{\al}^{\;2}}{2\big(e^{- c\la_L^{\al-}} - e^{-c
\la_L^{+}}\big)}.
\]
\el

\BPL \ref{Basic action estimate}. We have the identity
\[
E_c(U_L)=\; - \int_{- \infty}^{\la_L^{0-}} W^-(U_L) e^{cx}dx
\;+\;\int^{\infty}_{\la_L^{0-}} W^+(U_L) e^{cx}dx \;+\;\frac{1}{2}
\int_{\R}\big| (U_L)_x \big|^2 e^{cx} dx.
\]
We estimate each term separately, recalling that $W(U_L) \geq \al$
on $[\la_L^{\al-},\la_L^{+}]$ and $W^-(U_L) \leq W^-(a^-)$:
\begin{align*}
\int_{- \infty}^{\la_L^{0-}} W^-(U_L) e^{cx}dx \; & \leq \; W^-(a^-)
\int_{- \infty}^{\la_L^{0-}} e^{cx}dx \\
& =\; \frac{W^-(a^-)}{c}e^{c \la_L^{0-}},
\end{align*}
\begin{align*}
\int^{\infty}_{\la_L^{0-}} W^+(U_L) e^{cx}dx \; & \geq \;
\int^{\la_L^{+}}_{\la_L^{\al-}} W^+(U_L) e^{cx}dx \\
&\geq \; \al \int^{\la_L^{+}}_{\la_L^{\al-}} e^{cx}dx \; = \;
\frac{\al}{c} \Big[e^{c \la_L^{+}} - e^{c \la_L^{\al-}} \Big],
\end{align*}
\begin{align*}
d_{\al} \;  \leq \; \big| U_L(\la_L^{\al-}) - U_L(\la_L^{+}) \big|
\; & \leq \; \int_{\la_L^{\al-}}^{\la_L^{+}} \big| (U_L)_x\big|
dx\\
& \leq \; \Bigg(\int_{\la_L^{\al-}}^{\la_L^{+}} e^{-cx} dx
\Bigg)^{\frac{1}{2}} \Bigg(\int_{\la_L^{\al-}}^{\la_L^{+}} \big|
(U_L)_x \big|^2 e^{cx} dx \Bigg)^{\frac{1}{2}}.
\end{align*}
Hence, we have
\[
d_{\al}^{\;2}\; \leq \; \Bigg(\frac{e^{- c\la_L^{\al-}} - e^{-
c\la_L^{+}}}{c} \Bigg) \int_{\R} \big| (U_L)_x \big|^2 e^{cx} dx.
\]
Putting these bounds together, we obtain the desired estimate. \qed

\ms

\noindent \textbf{The speed of the travelling wave.} Thus far, all
the results were valid for an arbitrary $c>0$. It is easy to see
that the specific $c=c^*$ that guarantees existence should be very
special: by Proposition \ref{Approximate relation for c},
\begin{align*}
 \left(\big|(U_L)_x(+L^{+})\big|^2 -  \big|(U_L)_x(+
L^{-})\big|^2\right) &+  \left(\big|(U_L)_x(-L^{+})\big|^2
-\big|(U_L)_x( -L^{-})\big|^2\right)& \\
  +  \; 2  W^-(a^-) & = \;  2c \int _{\R} \big|(U_L)_x\big|^2 dx \\
& \geq\;  2c \int _{-L}^L \big|(U_L)_x\big|^2 dx \\
&   \geq \;  \frac{c \big|U_L(+L) - U_L(-L)\big|^2}{L} \\
& \geq \;  \frac{c}{L}\Big(|a^+ -a^-|-2r_0\Big)^2,
\end{align*}
which shows that if $c \larrow + \infty$ we can not achieve the
smooth matching of piecewise solutions at any $L < \infty$. On the
other hand, by Corollary \ref{The action equals derivatives' jump
discontinuities} and the a priori bound (\ref{Affine function
bound}), we have
\begin{align*}
{e^{+cL}} \left(\big|(U_L)_x(+L^{+}) \big|^2 -\big|(U_L)_x(+
L^{-})\big|^2\right) &+ {e^{-cL}}\left(\big|(U_L)_x(-L^{+})\big|^2
-\big|(U_L)_x( -L^{-}) \big|^2 \right) \\ & =\;2 c  E_c(U_L)\\
& \leq \;  2 c E_c(U_{aff})\\
& \leq \; -2e^{-c} W^-(a^-) + 2 c  e^{c} \big(E^+_0(U_{aff})\big),
\end{align*}
which shows that derivatives can not match if $c \larrow 0^+$. The
desired $c=c^*$ is the specific value, at which, for sufficiently
large $L > L^* \geq 1$, $E_c(U_L)=0$. This behavior of $E_{c>0}$ is
not present in its $E_{c=0}$ counterpart (\cite{A-F}, \cite{A-Be-C})
but it is \emph{plausible}: $\eq$ \emph{is} translation invariant
while (\ref{action}) \emph{is not}. Translates $U(\cdot - \de)$,
$\de \neq 0 $ of solutions occur as minimizers to a rescaled
$e^{c\de}E_c$, but both waves have the \emph{same} \emph{action}
only if $E_c(U(\cdot - \de)) = E_c(U)=0$.

\begin{remark} \label{Action non-increasing in L}
Note that for fixed $c>0$, the function $L \lmapsto E_c(U_L) :
[1,\infty) \larrow \big(-\infty, E_c(U_{aff})\big]$ is
\emph{non-increasing} in $L$: as $L$ increases, $\mathcal{X}_L$
increases ($L < L'$ implies $\mathcal{X}_L \subset
\mathcal{X}_{L'}$) and $E_c(U_L)$ decreases (see Sec.\ 2 for
definitions).
\end{remark}

\noi The next two estimates are key ingredients and will allow
determine of the speed and establish existence. The full strength of
\eqref{h*} is employed to show that $U_L$ \emph{can not get trapped
for infinite time} inside $\Ca^-$, after exiting the ball
$\mathbb{B}(a^-,r_0)$. We set
\[
R_{\max}^{\al}:= \underset{u \in
\partial \Ca^-}{\max}\ \big|u-a^- \big|.
\]

\bl \label{Uniform time bound below W=a} If $W$ satisfies ($h^*$),
there exists a $w^*
>0$ such that if $\al \in [0,\bar{\al}_0]$,
\begin{align}
\la_L^{\al-} - \la_L^{-} & \leq \  \frac{1}{w^*}\left\{c
R_{\max}^{\al} + \Big[\big(c R_{\max}^{\al} \big)^2 + 2 w^* \big|
R_{\max}^{\al} - r_0 \big| \Big]^{\frac{1}{2}}\right\}\nonumber\\
 & =: \ \La_{\al,-}.\label{Lambda Bound below W=a}
\end{align}
 As $w^*$ we may
take
\[
w^*\;:=\; \underset{| \xi | =1}{\underset{r_0 \leq r \leq
R_{\max}^{\al}}{\min} } \left[ \frac{d}{dt}\Big|_{t=r} W(a^{-} + t
\xi) \right].
\]
 \el

\BPL \ref{Uniform time bound below W=a}. Writing $\eq$ in polar form
$U_L =  a^- + \rho^-_L n^-_L$, we get (\ref{Equation in polar}).
Employing $(2)$ of \eqref{h*} on $[\la_L^{-},\la_L^{\al-}] \subseteq
[-L,L]$, we estimate
\begin{align*}
(\rho^-_L)_{xx}+ c (\rho^-_L)_x \; & \geq \; \nabla W(a^{-} +
\rho^-_L n) \cdot n^-_L \\
& = \; \frac{d}{dt}\Big|_{t=\rho^-_L} W(a^{-} + t
n^-_L)\\
& \geq \; \underset{| \xi | =1}{\underset{r_0 \leq r \leq
R_{\max}^{\al}}{\min}} \left[ \frac{d}{dt}\Big|_{t=r} W(a^{-} + t
\xi) \right] \;=:\; w^* \;>\; 0.
\end{align*}
Integrating once $(\rho^-_L)_{xx}+ c (\rho^-_L)_x \; \geq \; w^*$ on
$[\la_L^{-},x]$, $x\leq \la^{\al-}_L$ we get
\[
(\rho^-_L)_x + c \rho^-_L \; \geq \; w^*(x - \la_L^{-}) \;+\;
\Big\{c \rho^-_L(\la_L^{-}) \;+\;
(\rho^-_L)_x\big((\la_L^{-})^+\big)\Big\}.
\]
By Proposition \ref{Global control on minimizers}, we have $\{
\;.\;\} \geq 0$. By a further integration,
\[
\int_{\la_L^{-}}^x (\rho^-_L)_z(z)dz \;+\; c \int_{\la_L^{-}}^x
(\rho^-_L)(z)dz \;\geq \; w^* \int_{\la_L^{-}}^x (z - \la_L^{-})dz.
\]
Set $x := \la_L^{\al-}$. We utilize the a priori bound $
\big\|\rho^-_L \big\|_{L^{\infty}[\la_L^{-},\la_L^{\al-}]} \leq
R_{\max}^{\al}$ and that the right term equals $\frac{w^*}{2}\big[
\la_L^{\al-} - \la_L^{-} \big]^2$ to obtain
\[
\frac{\big| R_{\max}^{\al} - r_0 \big|}{\la_L^{\al-} - \la_L^{-}}
\;+\; \frac{c}{\la_L^{\al-} - \la_L^{-}} \Bigg( \int_{\la_L^{-}}^x
(\rho^-_L)(z)dz \Bigg) \;\geq \; \frac{w^*}{2}\big[ \la_L^{\al-} -
\la_L^{-} \big],
\]
which gives the desired inequality. Setting $\la_L^{\al-} -
\la_L^{-}=:x$ and comparing with the solutions of the parabola
$\frac{w^*}{2}x^2 - (c R_{\max}^{\al}) x - \big| R_{\max}^{\al} -
r_0 \big| \leq 0$ we obtain
\[
\frac{w^*}{2}\big[ \la_L^{\al-} - \la_L^{-} \big] \; \leq \;
\frac{\big| R_{\max}^{\al} - r_0 \big|}{\la_L^{\al-} -
\la_L^{-}}\;+\; c R_{\max}^{\al}
\]
which clearly implies \eqref{Lambda Bound below W=a}. \qed

\bl \label{Uniform time bound above W=a} For all $\al \in (0,
\bar{\al}_0]$, we have the implication:
\begin{equation} \label{Lambda Bound above W=a}
E_c(U_L) \leq  0 \ \ \ \Longrightarrow \ \ \ \la_L^{+} -
\la_L^{\al-} \; \leq\;
\frac{1}{c}\ln\left(1+\frac{W^-(a^-)}{\al}\right) \; =: \;
\La_{\al,+}.
\end{equation}
\el

\BPL \ref{Uniform time bound above W=a}. Follows directly from the
estimate of Lemma \ref{Basic action estimate}:
\begin{align*}
0\; & \geq \;E_c(U_L) \\
 & \geq \; e^{c \la_L^{\al-}} \Bigg\{-
\frac{W^-(a^-)}{c} \;+\; \frac{\al}{c} \Big(e^{c (\la_L^{+}-
\la_L^{\al-})} - 1 \Big) \;+\; \frac{c\; d_{\al}^{\;2}}{2\big(1 -
e^{- c (\la_L^{+}- \la_L^{\al-})}\big)} \Bigg\}\\
&  \geq \; \frac{ e^{c \la_L^{\al-}}\al}{c} \Big\{-\Big(
\frac{W^-(a^-)}{\al}\ + \ 1\Big)\;+\; e^{c (\la_L^{+}-
\la_L^{\al-})}\Big\}. \qed
\end{align*}

\begin{corollary} \label{bounds} The length of the time interval $ \big[\la^-_{L} ,\la_L^+\big]$ for which
the graph of $U_L$ remains between the constraint cylinders is $L$ -
uniformly bounded as long as $E_c(U_L)\leq 0$.
\end{corollary}

\BPCO \ref{bounds}. By Lemmas \ref{Uniform time bound below W=a} and
\ref{Uniform time bound above W=a}, we have
 \begin{align}
\la^+_{L} - \la^-_{L} \;
             & = \; \big( \la^{+}_{L} -
                          \la^{\al-}_{L}\big) \;+\;
                           \big(\la^{\al-}_{L} - \la^-_{L} \big)\nonumber\\
              & \leq \; \Lambda_{\al,+} \;+\; \Lambda_{\al,-}  \label{Lambda bound}\\
              & =: \; \Lambda \ <\ \infty,\nonumber
 \end{align}
provided that $E_c(U_L)\leq0$. This proves the bound. \qed

\begin{proposition} \textbf{(Determination of the speed of the travelling wave)}
\label{Determination of $c=c^*$} There exist  $c^*
>0 $ and $L^* \geq 1$ such that, for all $L \geq L^*$,
\[
E_{c^*}(U_L) \ =\  \underset{\mathcal{X}_L }{\inf}\ [E_{c^*}] \; =
\; 0.
\]
\end{proposition}

\noi The proof consists of several lemmas.

 \bl \label{c1} For any $L\geq 1$ and any $V \in \mathcal{X}_L$,
 both fixed, the function $c \lmapsto E_c(V)$ is continuous on $F:=\{c>0:|E_c(V)|<\infty\}$.
  \el

 \BPL \ref{c1}. Let $c_m \larrow c_\infty >0$ as $m\larrow \infty$.
 Since $V \in \mathcal{X}_L$, we have $W(V)=W^+(V)\geq 0$ on
 $[L,\infty)$ and as a result, for any $c\in F$,
 \begin{align*}
0\ & \leq \ \int_L^\infty
             \left(\frac{1}{2}|V_x|^2+W(V)\right)e^{cx}dx\\
   & = \ E_c(V)\ - \ \int_{-\infty}^L
             \left(\frac{1}{2}|V_x|^2+W(V)\right)e^{cx}dx\\
   & \leq \ E_c(V) \ + \ \sup_{(-\infty,L]}|W(V)|\int_{-\infty}^L
     e^{cx}dx\\
   & < \ \infty.
\end{align*}
Hence, for $m$ large we have on $(L,+\infty)$ that
\[
\left|\left(\frac{1}{2}|V_x|^2+W(V)\right)e^{c_m Id}\right|\ \leq \
2\left(\frac{1}{2}|V_x|^2+W(V)\right)e^{c_\infty Id} \ \in
L^1(L,+\infty).
\]
Again for any $c\in F$, we have
 \begin{align*}
\int_{-\infty}^L   \left|\frac{1}{2}|V_x|^2+W(V)\right|e^{cx}dx \ &
             \leq \ \int_{-\infty}^L
             \left(\left\{\frac{1}{2}|V_x|^2+W(V)\right\}
+\             2|W(V)|\right)e^{cx}dx\\
   & \leq \ E_c(V) \ + \ 2 \sup_{(-\infty,L]}|W(V)|\int_{-\infty}^L
     e^{cx}dx\\
   & < \ \infty.
\end{align*}
Since $c_m \larrow c_\infty$ as $m\larrow \infty$, if we choose $m$
large enough such that $c_m \leq \frac{3}{2}c_\infty$, we have
$e^{c_m x} \leq e^{c_\infty L} e^{\frac{c_\infty}{2} x}$ for all
$x\leq L$. Hence, for $m$ large we have on $(-\infty,L)$ that
\[
\left|\left(\frac{1}{2}|V_x|^2+W(V)\right)e^{c_m Id}\right|\ \leq \
e^{c_\infty L} \left(\frac{1}{2}|V_x|^2+W(V)\right)
e^{\frac{c_\infty}{2} Id} \ \in L^1(-\infty,L).
\]
By the pointwise convergence
$\left(\frac{1}{2}|V_x|^2+W(V)\right)e^{c_m Id} \larrow
\left(\frac{1}{2}|V_x|^2+W(V)\right)e^{c_\infty Id}$ as $m\larrow
\infty$, the lemma follows by application of the Dominated
convergence theorem on $(-\infty,L)$ and $(L,+\infty)$ separately.
 \qed

\ms

\noi Recall that $U_L$ has so far always denoted the minimizer of
$E_c$ into $\mathcal{X}_L$ for fixed $c$. We will temporarily denote
the dependence of $U_L$ on $c$ explicitly by $U_{L,c}$. Following an
idea of Heinze \cite{Hei}, we introduce the following set
 \be \label{SET}
C\ := \ \Big\{c>0 \ \big| \ \exists \ L\geq 1 \ : \ E_c(U_{L,c})<0
\Big\}.
 \ee
  \bl
 \label{L1} The set \eqref{SET} is open, non-empty and $\sup\ C
\leq \sqrt{2W^-(a^-)}{d_0}^{-1}$.
 \el

\BPL \ref{L1}. By observing that $C$ equals the set
\[
\Big\{c>0\ \big| \ \exists \ L\geq 1 \ \& \ \exists \ V  \in
\mathcal{X}_L : E_c(V)<0 \Big\},
\]
Lemma \ref{c1} implies that $C$ is open. By the bound \eqref{Affine
function bound} on $U_{\textrm{aff}} \in \bigcap_{L\geq 1}
\mathcal{X}_L$, we have $f(c) \ \geq \ E_c(U_{\textrm{aff}})$, where
\[
f(c):=e^{-c}\left(-\frac{1}{c}W^-(a^-) +
e^{2c}E^+_0(U_{\textrm{aff}})\right).
\]
Moreover, the equation $f(c)=0$ has a unique solution $c_0>0$ since
$f$ changes sign and $f'>0$ on $(0,\infty)$. Hence, $(0,c_0) \sub C
\neq \emptyset$. Moreover, by Lemma \ref{Basic action estimate}, for
$c \in C$ fixed, we have
\[
0\;> E_c(V) \ \geq \;E_c(U_L) \; \geq \; e^{c \la_L^{\al-}} \Bigg[-
\frac{W^-(a^-)}{c} \;+\; \frac{c\; d_{\al}^{\;2}}{2\big(1 - e^{- c
(\la_L^{+}- \la_L^{\al-})}\big)} \Bigg].
\]
which implies that $0 \geq c^2 d_\al^2 - 2 W^-(a^-)$. Letting $\al
\larrow 0^+$, we finally obtain
\[
0<c_0 \ \leq \ \sup C \ \leq \ \sqrt{2W^-(a^-)}{d_0}^{-1}.  \qed
\]
\medskip

 \bl \label{c2} Suppose that $L\geq 1$ is fixed and we have a sequence $C \ni c_m
 \larrow c_\infty$  as $m\larrow \infty$, $c_\infty >0$. Then, there exists a subsequence
  $c_{m,k}\larrow c_\infty$ along which
 \[
E_{c_{m,k}}(U_{L,c_{m,k}})\ \larrow  \ E_{c_\infty}(U_{L,c_\infty}),
\text{ \ \ as \ \ $k \larrow \infty$}.
 \]
 \el

\BPL \ref{c2}. Fix $\e >0$ and choose $V \in \mathcal{X}_L$ such
that $E_{c_\infty}(V)-\e \leq E_{c_\infty}(U_{L,c_\infty}) \leq
E_{c_\infty}(V)$. Since $c_m \larrow c_\infty$, by Lemma \ref{c1},
we can choose $m(\e) \in N$ large such that $|E_{c_\infty}(V)-
E_{c_m}(V)|\leq \e$, for all $m\geq m(\e)$. Thus,
 \begin{align*}
E_{c_{m}}(U_{L,c_{m}})\ & \leq \ E_{c_{m}}(V)\\
                        & \leq \ E_{c_\infty}(V) \ + \e \\
                        & \leq \ E_{c_\infty}(U_{L,c_\infty})  \ +
                        2\e,
\end{align*}
which implies
 \be \label{c3}
\limsup_{m \rightarrow \infty} \ E_{c_{m}}(U_{L,c_{m}})\  \leq \
E_{c_\infty}(U_{L,c_\infty}).
 \ee
By arguing as in the proof of Theorem \ref{Constrained Problem},
there exists a subsequence $c_{m,k}\larrow
 c_\infty$ along which $U_{L,c_{m,k}} \larrow \overline{U}$ in $[C^0_{\text{loc}}(\R)]^N$ and $U_{L,c_{m,k}}
 \lharpoonup \overline{U}$ weakly in $[H^1_{\text{loc}}(\R)]^N$, as $k \larrow
 \infty$. By weak LSC of the $L^2$ norm, we have
\[
\liminf_{k \rightarrow \infty} \ \frac{1}{2}\int_\R
|(U_{L,c_{m,k}})_x|^2 e^{c_{m,k}x}dx \ \geq \  \frac{1}{2}\int_\R
|(U_{L,c_{\infty}})_x|^2 e^{c_{\infty}x}dx.
\]
For $k$ large, we have the lower bound
\[
W(U_{L,c_{m,k}})e^{c_{m,k} Id} \geq -\big(e^{c_\infty
L}W^-(a^-)\big)e^{\frac{c_\infty}{2} Id} \chi_{(-\infty,L]}
\]
which is an $L^1(\R)$ function. Hence, the Fatou lemma implies
\[
\liminf_{k \rightarrow \infty} \ \int_\R W(U_{L,c_{m,k}})
e^{c_{m,k}x}dx \ \geq \ \int_\R W(U_{L,c_{\infty}})
e^{c_{\infty}x}dx.
\]
We conclude that
 \be \label{c4}
\liminf_{k \rightarrow \infty} \ E_{c_{m,k}}(U_{L,c_{m,k}})\  \geq \
E_{c_\infty}(\overline{U})\  \geq \ E_{c_\infty}(U_{L,c_\infty}).
 \ee
Putting \eqref{c3} and \eqref{c4} together, the proof follows. \qed

 \bl  \label{L2} If $c^*:=\sup\ C$, then $E_{c^*}(U_{L,c^*})=0$ for all $L\geq \La$.
 \el

\BPL \ref{L2}. By \eqref{SET}, there exists a sequence $C \ni c_m
\larrow c^*$ as $m \larrow \infty$ such that $E_{c_{m}}
(U_{L_m,c_{m}}) <0$. By the negativity of the action we may employ
the bound \eqref{Lambda bound} to obtain
\[
\la_{L_m}^+ - \la_{L_m}^- \ \leq  \ \La
\]
which is uniform in $m \in \N$. Moreover, since
$E_{c_{m}}(U_{L_m,c_{m}})<0$, we necessarily have $\la_{L_m}^+=L_m$,
since otherwise a translation to the right would contradict
minimality of $U_{L_m,c_{m}}$. By observing that the translate
$U_{L_m,c_{m}}(\cdot + L_m)$ is in $\mathcal{X}_\La$, we have
 \begin{align*}
E_{c_{m}}(U_{\La ,c_{m}})
       \ & \leq \ E_{c_{m}}(U_{L_m,c_{m}}(\cdot + L_m))\\
       & = \ e^{-c_mL_m}E_{c_{m}}(U_{L_m,c_{m}})\\
       & <\ 0.
\end{align*}
By Lemma \ref{c2}, the passage to the limit as $m \larrow \infty$
(along a subsequence if necessary) implies
\begin{align*}
E_{c^*}(U_{\La ,c^*})
   \ & = \ \lim_{m\rightarrow \infty} \ E_{c_{m}}(U_{\La ,c_{m}})\\
     &  \leq\ 0.
\end{align*}
Since $c^*=\sup\ C$ and $C$ is open, $c^* \notin C$ and as a result
$E_{c^*}(U_{\La ,c^*}) \geq 0$. By Remark \ref{Action non-increasing
in L} and \eqref{SET}, we conclude that $E_{c^*}(U_{L,c^*})=0$ for
all $L\geq \La$. \qed

\medskip

\BPP \ref{Determination of $c=c^*$}. By putting Lemmas \ref{c1},
\ref{L1}, \ref{c2} and \ref{L2} together, the proof of Proposition
\ref{Determination of $c=c^*$} follows with $c^*=\sup  C$,
$L^*=\La$. \qed

\ms

\noi Proposition \ref{Determination of $c=c^*$} provides a $c^*$ for
which $E_{c^*}(U_L)=0$ for large $L$ and this is sufficient for
existence. However, $c^*$ is the unique possible speed of minimizing
travelling waves\footnote{this fact together with a sketch of its
proof has been kindly pointed out by the referee.}:

\begin{proposition} \label{Uniqueness of the speed} \textbf{(Uniqueness of the speed)} Assume that
a minimizing solution $(U,c)$ to \eqref{problem} exists. Then, there
exists precisely one constant $c_*$ such that $(U,c_*)$ solves
\eqref{problem}.
\end{proposition}

\begin{corollary} Since minimizers of \eqref{action} have
vanishing action, we have $c_*=c^*$. Hence, Proposition
\ref{Determination of $c=c^*$} provides the unique constant for
which $E_{c^*}(U)=0$.
\end{corollary}

\BPP \ref{Uniqueness of the speed}. Let $(U_1,c^*_1)$, $(U_2,c^*_2)$
be two solutions of \eqref{problem} with $0< c^*_1 < c^*_2$ and
possibly $U_1 =U_2$. The differential form of the formula in Lemma
\ref{1st energy integral} is
\[
\frac{|U_x|^2}{2}+W(U)\ =  \ e^{-cx}\left( \frac{e^{cx}}{c}\left[
W(U) - \frac{|U_x|^2}{2} \right]\right)_x.
\]
We set $c:=c^*_2$, $U:=U^*_2$, multiply by $e^{c^*_1x}$ and
integrate by parts the right hand side to obtain
\begin{align*}
\int_{-t}^te^{c^*_1 x}\left(\frac{|(U_2)_x|^2}{2} + W(U_2)\right)dx
\ = \ & \left( \frac{e^{c^*_1 x}}{c^*_2}\left[W(U_2) -
\frac{|(U^*_2)_x|^2}{2}\right] \right)\Bigg|_{-t}^t \\
& -(c^*_1 - c^*_2)\int_{-t}^t  \frac{e^{c^*_1 x}}{c^*_2} \left[
W(U_2) - \frac{|(U_2)_x|^2}{2} \right]dx.
\end{align*}
We rewrite this identity as

\begin{align*}
\left({e^{c^*_1 x}}\left[W(U_2) - \frac{|(U^*_2)_x|^2}{2}\right]
\right)\Big|_{-t}^t
\ =\ \ & c^*_2 \int_{-t}^te^{c^*_1 x}\left(\frac{|(U_2)_x|^2}{2} + W(U_2)\right)dx \\
& +\ c^*_1 \int_{-t}^t  {e^{c^*_1 x}} \left[ W(U_2) -
\frac{|(U_2)_x|^2}{2} \right]dx\\
& -\ c^*_2 \int_{-t}^t  {e^{c^*_1 x}} \left[ W(U_2) -
\frac{|(U_2)_x|^2}{2} \right]dx\\
=\ \ & c^*_2 \int_{-t}^te^{c^*_1x}{|(U_2)_x|^2}dx\\
& -\ c^*_1 \int_{-t}^t  {e^{c^*_1 x}} \frac{|(U_2)_x|^2}{2} dx \ +\
c^*_1 \int_{-t}^t  {e^{c^*_1 x}} W(U_2) dx\\
= \ \ & (c^*_2 - c^*_1)\int_{-t}^te^{c^*_1 x} |(U_2)_x|^2dx \ + \
c^*_1 E_{c^*_1}\big(U_2 ,(-t,t)\big).
\end{align*}
Hence, we have the identity
\[
c^*_1 E_{c^*_1}\big(U_2 ,(-t,t)\big)\ =  \ (c^*_1 -
c^*_2)\int_{-t}^t|(U_2)_x|^2e^{c^*_1 x}dx \ + \ \left( e^{c^*_1
x}\left[W(U_2) - \frac{|(U_2)_x|^2}{2}\right] \right)\Big|_{-t}^t.
\]
By Proposition \eqref{Asymptotic behavior of constraint minimizers},
$(U_2)_x \larrow 0$ as $t \rightarrow \pm \infty$ up to sequences.
Since $E_{c^*_2}(U_2)=0$, we have
\begin{align*}
\int_\R \left\{\frac{|(U_2)_x|^2}{2}+ W^+(U_2) \right\}e^{c^*_2 x}dx
\
&= \ \int_\R W^-(U_2) e^{c^*_2 x}dx \\
& \leq \ W^-(a^-)\frac{e^{c^*_2 L_2^*}}{c^*_2}\\
& <\ \infty.
\end{align*}
where $L_2^*$ is a large constant as in Proposition
\ref{Determination of $c=c^*$}. Hence, since $c^*_1 <c^*_2$ we may
let $t \rightarrow \infty$ to obtain
\[
c^*_1 E_{c^*_1}(U_2)\ =  \ (c^*_1 - c^*_2)\int_\R|(U_2)_x|^2e^{c^*_1
x}dx \ < 0.
\]
But this contradicts that $c^*_1 E_{c^*_1}(U_2)\geq 0$. \qed

\noi We therefore in the remaining assume that $c=c^*$, the unique
speed provided by Proposition \ref{Determination of $c=c^*$}.

\medskip

\noindent \textbf{A Variational characterization of minimizing
travelling waves.} Summarizing, solutions $(U,c)$ to the system of
equations
\[
E_c(U) = \inf\ \Big\{E_c(V) :\ V \in [H_{\textrm{loc}}^1(\R)]^N,\
V(\pm \infty )= a^{\pm} \Big\} \; , \;\;\;\; \ E_c(U) = 0,\\
\]
are heteroclinic travelling waves and solve the differential
equations
\[
\left\{\begin{array}{l}
  \eq\\
   U(\pm \infty) = a^{\pm}.\\
\end{array}\right.
\]
Both the weight $e^{c \text{Id}}$ of \eqref{action} and its
minimizer are unknown. The first equation of the system involves the
minimization problem for $E_c$ in the class $\{E_c\;|\;c>0\}$ and,
the second one selects $c=c^*$ so that the minimum zero.

\section{Removing the Constraints.}

In this section we prove existence of solution to problem
\eqref{problem}.

\bt \textbf{\emph{(Existence)}} \label{Existence result} Assume the
potential $W$ satisfies \eqref{h}, \eqref{h*}. Then, there exists a
travelling wave solution $(U,c) \in [C^2(\R)]^N \times (0,+ \infty)$
to
\[
\left\{\begin{array}{l}
  \eq\\
   U(\pm \infty) = a^{\pm}.\\
\end{array}\right.
\]
The speed $c$ equals the constant $c^*$ in Proposition
\ref{Determination of $c=c^*$} which is unique. In particular,
$E_{c^*}(U)=0$. \et

\BPT \ref{Existence result}. By Proposition \ref{Determination of
$c=c^*$}, we have $E_{c^*}(U_{L,c^*})=0$, for all $L\geq L^*$. By
Corollary \ref{bounds}, if we choose $L > \Lambda$ we obtain a
minimizer $U := U_L$ of $E_c$ with $c=c^*$ for which $E_c(U)=0$.
Thus either $U$ or a translate $U(\cdot - \de)$ (with necessarily
the same action) does not realize the constraint, solving
(\ref{problem}) on $\R$. The proof is complete. \qed

\ms

\begin{corollary} The speed $c^*$ has the variational characterization\footnote{Analogous characterizations have been
obtained in \cite{H-P-S} and \cite{He} for other travelling wave
problems.}
\[
c^*\ = \ \sup_{c>0}\ \left\{c\ \Big| \inf_{V \in \mathcal{X}}\
E_c(V)<0 \right\},
\]
where $\mathcal{X}:=\big\{V \in [H_{\textrm{loc}}^1(\R)]^N: V(\pm
\infty )= a^{\pm} \big\}$.
\end{corollary}

\noindent We now derive a priori bounds on $c^*$. We take $t>0$ and
consider the affine $[W^{1,\infty}_{\textrm{loc}}(\R)]^N$ function
\begin{equation}
\label{Affine function with t} U^t_{aff}(x)\;:=\; a^- \chi_{(-
\infty , -t)} \;+\; \left( \frac{t - x}{2t} a^- + \frac{t + x}{2t}
a^+\right)\ \chi_{[-t , t]} \;+\;a^+ \chi_{(t, \infty)}.
\end{equation}

\begin{proposition} \textbf{(A priori bounds on $c^*$)} \label{A priori
bounds on the speed c} There exist $ 0 < c_{\min} < c_{\max} <
\infty$ depending only on $W$, such that
\[
c_{\min} \; \leq \; c^* \; \leq \; c_{\max}.
\]
Moreover, if $d_0 := {\lim}_{\al \rightarrow 0^+} d_{\al}$, then
\[
c_{\max}\;=\;\frac{\sqrt{2 W^-(a^-)}}{d_0},
\]
\[
 c_{\min}\;=\;
\underset{t >0}{\sup} \Bigg[ \frac{W^-(a^-)}{e^{2t c_{\max}}}
\Bigg(\frac{1}{2}\Big\{\frac{\big|a^+ -a^-\big|}{2t}\Big\}^2 +
\int_{-t}^{t} W^+ \Big(\frac{t - x}{2t} a^- + \frac{t + x}{2t} a^+
\Big) dx \Bigg)^{-1}\Bigg].
\]
\end{proposition}

\BPP \ref{A priori bounds on the speed c}. The upper bound follows
by Lemmas \ref{L1} and \ref{L2}. For the lower bound, we utilize
(\ref{Affine function with t}) and take as we can $t=L$. This gives
as in (\ref{Affine function bound}) that the inequality $0 =
E_c(U_t) \leq E_c(U^t_{aff})$ implies
\[
0 \; \leq \; -e^{-ct} \frac{ W^-(a^-)}{c} \;+\; e^{ct}
\int_{-t}^{t}\Big\{ \frac{1}{2} \Big| \frac{a^+ - a^-}{2t} \Big|^2 +
W^+ \Big( \frac{t - x}{2t} a^- + \frac{t + x}{2t} a^+ \Big) \Big\}
dx.
\]
Hence, for all $t > 0$,
\[
c\; \geq \; \frac{W^-(a^-)}{e^{2 c t}} \Bigg( \int_{-t}^{t}\Big\{
\frac{1}{2} \Big| \frac{a^+ - a^-}{2t} \Big|^2 + W^+ \Big( \frac{t -
x}{2t} a^- + \frac{t + x}{2t} a^+ \Big) \Big\} dx \Bigg)^{-1}.
\]
Utilizing the upper bound and maximizing with respect to $t>0$, we
are done.  \qed

\section{Extensions.}

\noi Utilizing ideas related to those in \cite{A-F}, we relax
\eqref{h*} to a localized version. The new \eqref{h**} requires the
existence of two convex components $\Ca^{\pm}$ of the sublevel set
$\big\{W \leq\al \big\}$, but only when $W$ is restricted in a large
convex $\Omega \subseteq \R^N$ without any restriction on
$W|_{\textrm{ext}(\Omega)}$. As a consequence, \eqref{h**} allows
for potentials with several other minima and/or unbounded values to
$- \infty$.
\[
\label{h**} \tag{$h^{**}$}
\begin{array}{l}\text{There exists a $C^2$
   convex closed set $\Omega \subseteq \R^N$
   which encloses the minima}\\
 \text{$a^\pm$ and satisfies (H3), such that \eqref{h*}
    holds for $W$ within $\Omega$. Moreover,}\\
    \text{The values of W on $\partial \Omega$ exceed those in the interior :
 if $u \in int(\Omega)$, then}\\
  \text{$W(u) < {\min}_{\partial \Omega} W$.}
\end{array}
\]

\begin{example}
\textbf{(N-d potentials satisfying \eqref{h*}, \eqref{h**})} (i) We
construct a deformation of the 2-well potential $W(u):=|u - a^+|^p
|u-a^-|^p$, $p\geq2$, $u \in \R^N$. We take $\e
> 0$ and set
\[
F_{\e}(u):= \big\{\e \exp\big((|u - a^-|^2 - \delta^2)^{-1}\big)
+1\big\} \chi_{\mathbb{B}(a^-, \delta)}\ + \ \chi_{\R^N \set
\mathbb{B}(a^-, \delta)},
\]
where $C:={\max}_{|u -a^-|=\delta}\{W(u)\}$ and define $W_{\e}(u):=
F_{\e}(u)W(u)-C(F_{\e}(u)-1)$. The potentials $W_{\e}$ satisfy our
assumptions and $W_{\e} - W \larrow 0$ in $C^2(\R^N)$, as $\e
\rightarrow 0^+$.

\noi (ii) \label{2D-Example} (G. Paschalides) The following
deformation of the 2-well planar potential
\[
W_C(u_1,u_2)\;:=\; \left\{\begin{array}{l}
  W(u_1,u_2)\;, \;\;\;\; \ \hspace{110pt} u_1 <0\;,\;u_2 \in \R,\\
  W(u_1,u_2) - C\big[6u_1^5 - 15u_1^4 +10u_1^3 \big]\;,
 \;\;\;\;\;\; 0 \leq u_1 \leq 1\;,\;u_2 \in \R,\\
   W(u_1,u_2) - C\;, \;\;\; \ \hspace{95pt} u_1 > 1\;,\;u_2 \in \R,\\
\end{array}\right.
\]
(with $a^{\pm} = (\pm 1, 0)$) satisfies the assumptions \eqref{h*},
\eqref{h**} for any $C>0$.
\end{example}

\begin{remark} Can monotonicity of \eqref{h*},
\eqref{h**} be relaxed? In the Appendix we construct a class of
$W$'s which are monotone except for merely one critical point $a^0$
in $W^{-1}([W(a^-),0])$. This implies existence of a connection $a^+
- a^0$, different from $a^+ - a^-$, which generally obstructs
existence. Critical points at lower level attract, for $c>0$, the
flow of $\eq$ (see also Risler \cite{R}).
\end{remark}

\ms

\noindent \textbf{Extension of Theorem \ref{Existence result} under
the assumption \eqref{h**}.}  In this case we solve a related
problem for a modified ``better'' $\overline{W}$ and then show that
the solution we construct is also a solution of the original problem
as well. We modify $W$ to a new $\overline{W}$ by setting:
\[
\overline{W}\ :=\ W\chi_{\{W \geq {\min}_{\partial \Omega}W \}} \ +\
\big(2\ {\min}_{\partial \Omega}W - W\big)\chi_{\{W
<{\min}_{\partial \Omega}W \}}.
\]
This is the reflection the graph of $W$ with respect to the
hyperplane $\big\{w = {\min}_{\partial \Omega}W \big\}$ which maps
any parts of $Gr(W)$ lying into $\big\{W \leq {\min}_{\partial
\Omega}W \big\}$, to the opposite halfspace. $\overline{W}$ is
sufficiently coercive and Lemma \ref{l2} applied to $\Omega$ and to
$E_c$ provides an $[L^{\infty}(\R)]^N$ - bound for the minimizers,
showing that they are localized inside $\Omega$. Since
$\overline{W}$ satisfies \eqref{h*} inside $\Omega$, problem
(\ref{problem}) for $\overline{W}$ has a solution $U$ in
$[C^2(\R)]^N$. By construction $W \big|_{\Omega} \equiv \overline{W}
\big|_{\Omega}$, so $U$ solves (\ref{problem}) for $W$ as well.

\section{Appendix}

\noi \textbf{On the optimality of the assumptions.} We construct a
class of $W$'s for which there is a heteroclinic between a local
minimum $a^+$ with $W(a^+)=0$ and a critical point $a^0$ with
$0>W(a^0)>W(a^-)$, $a^-$ the global minimum. Hence, the existence of
additional solutions which may obstruct the existence of $a^+ - a^-$
connections can not be excluded without monotonicity as in
\eqref{h*}.
 \[
 \label{h1} \tag{h1}
 \text{We assume that $W \in C^2_{\textrm{loc}}(\R^N)$ and}\hspace{170pt}
\]
\begin{enumerate}
  \item $W$ has at least 3
critical points, $a^{\pm}$, $a^0$ with $a^{\pm}$ local minima, $a^0$
critical point and $W(a^+)=0 > W(a^0)> W(a^-)$.

\item For $2 \leq j \leq N$, $W_{u_j}(u_1,0,...,0)=0$ and
$[a^-,a^0]$, $[a^0,a^+]$ are on the $u_1$-axis.
\end{enumerate}

\noindent If $N=1$ and $a^0$ is a local minimum, then generally no
$a^+ - a^-$ connection exists  (\cite{F-McL}), depending on the
speeds $c_{-,0}$ and $c_{0,+}$ of the solutions $a^- - a^0$ and $a^0
- a^+$. For $N>1$, $(b)$ implies the existence of solutions
$U=(u,0,...,0)$ to $\eq$ for the slice
$\overline{W}(u):=W(u,0,...,0)$. Thus, we may only impose
assumptions on $\overline{W}$:
 \[  \label{h2} \tag{h2}\begin{array}{l}
 \text{We assume that \eqref{h**} holds, with the
exception that $\overline{W}$ is monotone on}\\
   \text{ $(a^-,a^0)$ $(a^0,a^+)$ separately, instead of $(a^-,a^+)$.}
\end{array}
\]

\noi\textbf{Proposition.} \emph{\label{bistable heteroclinics} If
$W$ satisfies \eqref{h1}, \eqref{h2}, there exists a solution $(U,c)
\in [C^2(\R)]^N\times (0,\infty)$ to
\[
\left\{\begin{array}{l}
  \eq\\
   U(+ \infty) = a^{+}\;,\;\;U(- \infty) = a^{0}.\\
\end{array}\right.
\]}
\BPP.We deform smoothly the slice $\overline{W}$ to a new
$\widehat{W}$ for which the \emph{nature of the critical point $a^0$
is changed}, being a global \emph{minimum} of ${\widehat{W}}$. Then,
the problem for $\widehat{W}$ can be tackled by the foregoing
theory, and, by a localization argument, the solution we construct
solves also the original problem. Let $F : (a^-,a^0) \larrow
(0,\infty)$ be the "half" of the standard bell function $F(u):=K
\exp \big((u - a^0)^{-1}(u - a^0 + 2 a^-)^{-1}\big)$, $K>0$ to be
chosen, and consider the following transformation
\[
 \widehat{W}(u)\;:=\;\left\{\begin{array}{l}
  \overline{W}(\Omega_2)\;,\hspace{85pt}\ \ \ \ \ u \geq \Omega_2\\
  \overline{W}(u)\;,\hspace{94pt}\ \ \ \  u \in [a^0, \Omega_2)\\
  -\big(F(u) \overline{W}(u) - 2 \overline{W}(a^-)\big)\;,
\;\;\;\;\;\;\;\; u \in (a^-, a^0)\\
  -\big(F(a^-) \overline{W}(u) - 2 \overline{W}(a^-)\big)\;,\;\;\;\;\;\; u \leq
a^-.\\
\end{array}\right.
\]
We choose $K>0$, such that $\;\widehat{W}(a^-) \geq
\overline{W}(\Omega_1)$. Assumptions \eqref{h1}, \eqref{h2} imply
that $\widehat{W}$ satisfies \eqref{h*}, giving an $a^+ \;-\; a^0$
heteroclinic which solves $u_{xx}-\widehat{W}'(u)=-c u_x$ (Theorem
\ref{Existence result}). Lemma \ref{l2} provides the
$L^{\infty}(\R)$ - bound
\[
a^- \; \leq \; u(x) \; \leq \;\Omega_2\;,\;\;\;\;\;\; \text{for all
} x \in \R.
\]
\[
\includegraphics[scale=0.16]{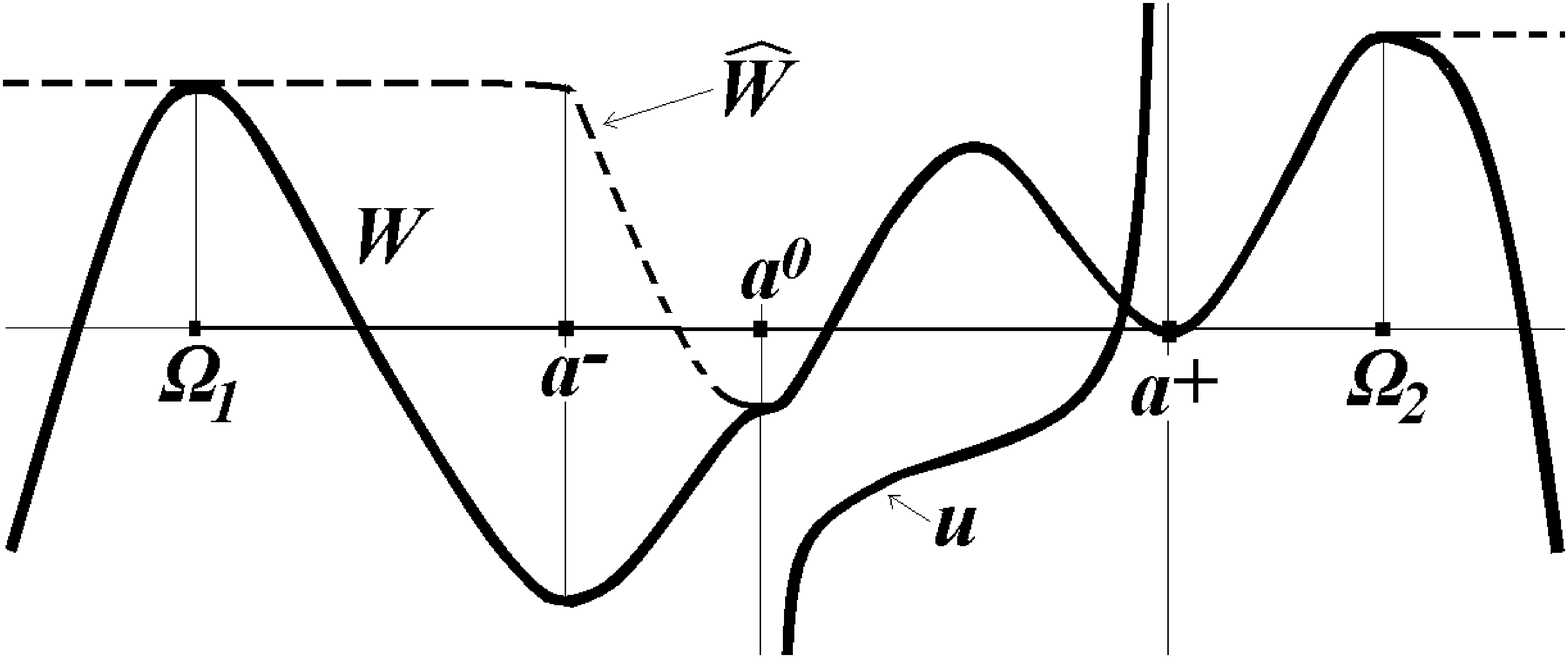}
\]
The function $u$ solves $u_{xx}-\overline{W}'(u)=-c u_x$ as well.
Indeed, it suffices to improve the bound on $u$ to $a^0 \leq u(x)
\leq a^+$, for all $x \in \R$. Since by construction
$\overline{W}\big|_{[a^0,a^+]} \equiv \widehat{W}\big|_{[a^0,a^+]}$.
Lemma \ref{l2} applied to (\ref{action}) for $\widehat{W}$ gives the
desired localization. \qed

\medskip
\noindent {\bf Acknowledgement.} {We thank Peter Bates, Vassilis
Papanicolaou and Achilleas Tertikas for their suggestions and their
interest in the present work. We also thank Gregory Paschalides for
the example \ref{2D-Example} and the numerical simulation. Special
thanks are due to the anonymous referee for his several suggestions
and valuable comments which improved the content as well as the
presentation of this paper. Finally, we wish to thank Hiroshi Matano
for the information he gave us on the status of the problem.}

\bibliographystyle{amsplain}

\end{document}